%% file: spacetime-MS.tex
\documentclass[final,onefignum,onetabnum]{siamart190516}
\ifpdf
\hypersetup{
  pdftitle={An entropy structure preserving space-time formulation for cross-diffusion systems},
  pdfauthor={M. Braukhoff, I. Perugia, P. Stocker}
}
\fi
\usepackage{a4wide}
\usepackage{hyperref}
\usepackage{graphicx}
\usepackage{cite}
\usepackage{color}
\usepackage{amsmath}
\usepackage{amssymb}
\usepackage{bm}
\usepackage{csvsimple}
\usepackage{siunitx} 
\usepackage{cleveref}
\usepackage{tabularx}
\usepackage{booktabs}
\usepackage{pgfplots} 
\pgfplotsset{compat=1.17}
\pgfplotsset{
discard if not/.style 2 args={
    x filter/.code={
        \edef\tempa{\thisrow{#1}}
        \edef\tempb{#2}
        \ifx\tempa\tempb
        \else
            
        \fi
    }
}
}
\pgfplotscreateplotcyclelist{paulcolors}{
teal,every mark/.append style={solid,fill=teal},mark=*,very thick,mark size=3pt\\
orange,every mark/.append style={solid,fill=orange},mark=square*,very thick,mark size=3pt\\
violet,every mark/.append style={solid,fill=violet},mark=diamond*,very thick,mark size=3pt\\
teal,densely dashed,every mark/.append style={solid,fill=teal},mark=*,very thick,mark size=3pt\\
orange,densely dashed,every mark/.append style={solid,fill=orange},mark=square*,very thick,mark size=3pt\\
violet,densely dashed,every mark/.append style={solid,fill=violet},mark=diamond*,very thick,mark size=3pt\\
teal,dash dot,every mark/.append style={solid,fill=teal},mark=*,very thick,mark size=3pt\\
}
\pgfplotscreateplotcyclelist{paulcolorsnm}{teal, very thick\\ orange,densely dashed, very thick\\ violet,dash dot, very thick\\ }

\input{pgfslopetriangle}

\renewcommand{\epsilon}{\varepsilon}
\newcommand{\WeightedSpace}{H^1_{\epsilon}}

\definecolor{pscol}{rgb}{0,0,0.7}

\definecolor{mbcol}{rgb}{0.7,0,0.7}

\definecolor{mbcolnew}{rgb}{0.9,0.1,0.2}

\definecolor{ipcol}{rgb}{0.9,0.0,0.9}

\newsiamremark{remark}{Remark}
\newsiamthm{example}{Example}

\title{An entropy structure preserving space-time {formulation} for
  cross-diffusion systems{: analysis and Galerkin discretization}}
\headers{A space-time method for cross-diffusion systems}{M. Braukhoff, I. Perugia, P. Stocker}

\author{Marcel Braukhoff\thanks{Mathematisches Institut, Heinrich-Heine-Universität Düsseldorf, Germany (\email{marcel.braukhoff@hhu.de}).} \and Ilaria Perugia\thanks{Faculty of Mathematics, University of Vienna, Austria (\email{ilaria.perugia@univie.ac.at}, \email{paul.stocker@univie.ac.at}).} \and Paul Stocker\footnotemark[2]}

\begin{document}  
\sloppy
\maketitle

\begin{abstract}
Cross-diffusion systems are systems of nonlinear parabolic partial
differential equations that are
used to describe dynamical processes in several application, including chemical concentrations and cell biology.
We present a space-time approach to the proof 
of existence of bounded weak solutions of cross-diffusion systems,
making use of the system entropy to examine long-term behavior and
to show that the solution is nonnegative, even when
a 
maximum principle is not available. 
This approach naturally gives rise to a novel space-time Galerkin
method for the numerical approximation of cross-diffusion systems
that conserves their entropy structure. 
We prove 
existence and convergence of the discrete solutions, and present
numerical results for the porous medium, the Fisher-KPP, and
the Maxwell-Stefan
problem. 
\end{abstract}
\begin{keywords}
Space-time Galerkin method, entropy method, strongly coupled parabolic systems, global-in-time existence, bounded weak solutions, space-time finite elements
\end{keywords}
\begin{AMS}
35K51, 
35K55, 
35Q92, 
65M60, 
41A10 
\end{AMS}

\section{Introduction}

In this paper we develop a new space-time approach to the celebrated
{\em boundedness by entropy method} by Ansgar J\"ungel
\cite{jungel2015boundedness}. For a textbook version see
\cite{jungel2016entropy}; see also \cite{equadiff,chen2015}.

Cross-diffusion systems are systems of nonlinear parabolic partial
differential equations that are
commonly used to describe dynamical processes appearing in modeling, for example, population dynamics, ion transport through nanopores, tumor growth models, and multicomponent gas mixtures.
The challenge in the analysis of these systems is that the diffusion matrix is not necessarily symmetric nor positive semi-definite, and thus no maximum principle is available.
Following \cite{jungel2015boundedness}, the remedy is to make use of the entropy structure of the system. 
Introducing the entropy function, a transformation of the solution, allows us to examine long-term behavior and show that the solution is nonnegative and bounded. 
Here, we present a space-time approach to the proof 
of existence of bounded weak solutions of cross-diffusion systems.
The main tool will be the method of compensated compactness, which is
a special technique {of} applying the classical div-curl lemma~\cite{TartarBook2009}.
The key difference to the existing literature is that we do not make
use of time-stepping, but instead consider time and space altogether.
This naturally leads to a novel space-time Galerkin method
for the numerical approximation of cross-diffusion systems.
The space-time approach entails test and trail spaces, as well as the mesh, where time is included as additional dimension. 
This provides an easy way to increase the approximation degree
simultaneously in space and time, and makes space-time $hp$-refinement
possible.
{In a schematic way, our overall approach consists of the following four
  steps:
  \begin{enumerate}
  \item space-time variational formulation,
  \item transformation to \emph{entropy} variables, 
  \item regularization with a space-time $H^1$ inner product,
  \item Galerkin discretization. 
  \end{enumerate}}

\medskip
Existing numerical schemes for cross-diffusion systems rely on time-stepping methods.
An entropy/energy conserving time-stepping algorithm for
thermomechanical problems
was developed in \cite{portillo2017energy} being of second order in time.
In \cite{kunstmann2018runge}, assuming existence of sufficient regular strong solutions on some time interval
$[0,T]$ of a scalar  diffusion equation, Runge-Kutta methods were studied using maximal regularity. 
Although maximal regularity also applies to a certain type of cross-diffusion systems \cite{pruss2016}, Runge-Kutta methods were only applied to very restrictive {classes}; an example (semi-discrete Runge-Kutta scheme) can be found in \cite{jungel2015entropy}. 
In \cite{junge2017fully}, an entropy diminishing/mass conserving fully
discrete variational formulation for a cross-diffusion system was
presented.
{An alternative discretization for cross-diffusion systems based on the
change to entropy-variables has been proposed in~\cite{Egger2019},
where a dissipation-preserving approximation by Galerkin methods in
space and discontinuous Galerkin methods in time has beed studied.}

\medskip
Maxwell-Stefan systems, see \cite{Stefan1871,Maxwell1866}, describe multicomponent diffusive fluxes in non-dilute solutions or gas mixtures, and are a prime example for the cross-diffusion systems considered here.
The first result on global solutions for the Maxwell-Stefan equations close to the equilibrium is given in \cite{GiMa1998}.
The global existence of solutions close to equilibrium and the large time convergence to this equilibrium can be found in \cite[Chapter 9]{Giovangigli2012}, \cite{GiMa1998asymptotic,HMPW2017}, and \cite[Chapter 12]{pruss2016}.
The proof of existence of local classical solutions to the Maxwell-Stefan equations can be found in \cite{Bothe2011}.  For a textbook on this topic, see \cite{pruss2016}.
The fact that the Maxwell-Stefan equations satisfy the assumptions made in this paper, see (H1)-(H3) below, is due to \cite{JS13}, where the entropy structure of the Maxwell-Stefan system was used to prove the existence of globally bounded weak solutions.
An entropy structure was also identified for a generalized Maxwell-Stefan system coupled to the Poisson equation in \cite{JuLe2019}, where the existence of global weak solutions was proven as well. 
The unconditional convergence to the unique equilibrium for given mass was shown in \cite{HMPW2017,MT2015} without reaction terms.
Those results were extended to also include reaction terms using mass-action kinetics in \cite{daus2018}, whenever a detailed balance equilibrium exists.
The heat equation can be recovered from the Maxwell-Stefan equation as a relaxation limit \cite{SaSo2018}.

{As of yet, numerical schemes for the Maxwell-Stefan equations commonly employ time-stepping.}
A finite differences approximation can be found in \cite{LeAn2009,LVM1992}. 
Fast solvers for explicit finite-difference schemes were studied in \cite{Geiser2015}.
A posteriori estimates for finite elements in the stationary case are given in \cite{carnes2008}. 
In \cite{peerenboom2011}, a mass conserving finite volume scheme was presented.
Existence of solutions for a mixed finite element scheme under some restrictions on the coefficients was proven in \cite{mcleod2014}.
The scheme of \cite{dieter2015} was proven to also conserve the $L^\infty$ bound by making use of a maximum principle. 
A scheme using finite elements in space and implicit Euler in time was used to approximate a Poisson-Maxwell-Stefan system in \cite{JuLe2019}. 
That scheme, which is based on a formulation in entropy variables, admits solutions that conserve the mass as well as the entropy structure. 
As a by product, the solution satisfy an $L^\infty$ bound.
Another scheme that is mass conserving and conserves the $L^\infty$ bounds of the solutions was presented in \cite{boudin2012}.

\medskip
{On simultaneous space-time finite element approaches for parabolic problems,
there is a rich literature on the linear case, focusing on the heat
equation; see, e.g. the recent overview~\cite[Ch.~7]{SpacetimeRICAM}.
We point out that due to the different orders of derivatives present,
conforming discretizations are typically based on a Petrov-Galerkin
approach; see e.g.~\cite{Aziz1989, Schwab2009,Andreev2013,Steinbach2015,Zank2021}.}
For nonlinear parabolic equations, the existing literature on space-time methods is much sparser.
The adaptive finite element scheme introduced in \cite{Eriksson1991}
for linear parabolic problems was extended in \cite{Eriksson1995} to the scalar version of the nonlinear reaction-diffusion equation treated in this paper.
A space-time discontinuous Galerkin method for scalar nonlinear
convection and diffusion was introduced in \cite{cesenek2012}.
A space-time method for nonlinear PDEs using adaptive wavelets was
introduced in \cite{MR2216616}.

\medskip
The structure of this paper is as follows.
In \cref{sec:generalsetting}, we state the problem and make the necessary assumptions for the existence of an entropy function.
In \cref{sec:STGM}, we present the space-time Galerkin method on a
regularized formulation of the problem in the entropy variable
unknown, and state our two main results in \cref{prop.nonlin} and \cref{prop2},
namely existence and convergence of discrete solutions, respectively.
Existence of discrete solutions is proven in \cref{sect:existence}.
The proof of convergence will be split into two parts, first showing
convergence with respect to mesh size in \cref{sec:limh}, then proving
convergence as the regularization parameter 
goes to zero in \cref{sec:limeps}.
In \cref{sect:existence_weak}, we are then able to prove existence of
a weak solution of the continuous problem.
Numerical tests for the porous medium, the Fisher-KPP, and the Maxwell-Stefan problem are presented in \cref{sec:numerics}.
All numerical results\footnote{The code is available online at \url{https://github.com/PaulSt/CrossDiff}} were obtained using the finite element software NGSolve, see \cite{joachim}. 
Additionally, in \cref{sec.MaxwellStefan}, we reformulate 
the Maxwell-Stefan system with implicitly given currents in terms
  of the concentrations, and test it numerically.
     
\section{General setting}\label{sec:generalsetting}
Let $\Omega\subset \mathbb R^d$ be a bounded domain, and $\rho_0\in
L^\infty(\Omega)^N$, $N\ge 1$, a vector-valued function. We
consider the following nonlinear reaction-diffusion system in the
vector-valued
unknown $\rho(t)=(\rho_1,\ldots,\rho_N)(\cdot,t):\Omega\to \mathbb
R^N$:
\begin{equation}
\label{eq.classical}
\begin{cases}
\partial_t \rho-\nabla \cdot(A(\rho)\nabla
\rho)=f(\rho)&\mbox{in }\Omega,\ t>0,\\
(A(\rho)\nabla\rho)\cdot\nu 
= 0&\mbox{on }\partial\Omega, \ t>0,\\
\rho(0)=\rho_0 &\mbox{in }\Omega.
\end{cases}
\end{equation}
Here, $A(\rho)\in \mathbb R^{N\times
N}$ is the diffusion matrix,
$f(\rho): \mathbb R^N\to \mathbb R^N$ represents the reactions,
and $\nu$ is the 
outward pointing unit normal vector at $\partial\Omega$;
moreover, for $1\le i\le N$,
\[
\left(\nabla\cdot(A(\rho)\nabla
\rho)\right)_i=\sum_{\mu=1}^d\sum_{j=1}^N\frac{\partial}{\partial
x_\mu}\left(A_{ij}(\rho)\frac{\partial \rho_j}{\partial
x_\mu}\right),\quad
\left((A(\rho)\nabla\rho)\cdot\nu\right)_i=\sum_{\mu=1}^d\sum_{j=1}^N
A_{ij}(\rho) \frac{\partial \rho_j}{\partial x_\mu}\nu_\mu.
\]

We make the following hypotheses, which are {slightly stronger assumptions compared to those} made by
A.~J\"ungel in~\cite{jungel2015boundedness}.
\begin{enumerate}
\item[(H1)]
$A\in C^{0}(\overline{\mathcal D};\mathbb R^{N\times N})$
and $f\in C^{0}(\overline{\mathcal D};\mathbb R^N)$, for a bounded
domain $\mathcal D\subset (0,\infty)^N$.
\item[(H2)] 
There exists a convex function $s\in C^2(\mathcal D,[0,\infty))\cap
C^0(\overline{\mathcal D})$, with
$s':\mathcal D\to\mathbb R^N$ invertible and $u:=(s')^{-1}\in
C^1(\mathbb R^N,\mathcal D)$,
such that the following two
conditions are satisfied:
\begin{enumerate}
\item[(H2a)]\label{con.A}
There exists a constant $\gamma>0$ such that
\begin{align*}\
z\cdot s''(\rho)A(\rho)z\geq \gamma |z|^2\qquad \forall z\in
\mathbb R^N,\,  
\rho\in\mathcal D.
\end{align*}
Note that $s''(\rho)$ is matrix-valued, with
$(s''(\rho))_{k\ell}=\frac{\partial}{\partial \rho_k}
(s'(\rho))_\ell=\frac{\partial^2}{\partial \rho_k \partial \rho_\ell}s(\rho)$.
\item[(H2b)] 
There exists a constant $C_f\geq0$ such that
\begin{align*}
f(\rho)\cdot s'(\rho)\leq C_f \qquad\forall  \rho\in\mathcal D.
\end{align*}
\end{enumerate}
\end{enumerate}
{Additionally, we make the following assumption on $\rho_0$:
\begin{enumerate}
\item[(H3)]  The initial datum $\rho_0:\Omega\to\overline{\mathcal D}$ is measurable. 
\end{enumerate}
}

A discussion on when it is possible to find a convex function $s$
such that $(H2)$ is satisfied for cross-diffusion equations can be
found in~\cite{chen2019} (see \cite[Lemma 22]{chen2019}).

Let $T>0$. {A weak formulation of \eqref{eq.classical} reads as
  follows:
  Find
$\rho\in L^2(0,T;H^1(\Omega)^N)\cap H^1(0,T;(H^1(\Omega)')^N)$ 
such that}
\begin{equation}
\label{eq.vweak}
\int_0^T\langle\phi ,\partial_t\rho\rangle dt
+ \sum_{i,j=1}^N\int_0^T\int_{\Omega}\nabla\phi_i \cdot A_{ij}(\rho)\nabla \rho_j dx dt=\int_0^T\int_{\Omega}\phi\cdot f(\rho) dxdt
\end{equation}
for all $\phi\in L^2(0,T;H^1(\Omega)^N)$, with
    $\rho(0)=\rho_0$,
where $\langle \cdot,\cdot\rangle$ denotes the duality product
between $H^1(\Omega)^N$ and $(H^1(\Omega)')^N$.

{By introducing the so-called entropy variable $w$, which satisfies
  $\rho=u(w)$, problem~\eqref{eq.classical}, as well as~\eqref{eq.vweak}, can be rewritten in terms in
  the unknown $w$.}
{
    \begin{remark}\label{rem:as}
          In~\cite{jungel2015boundedness},
          a more degenerate version of (H2a) is permitted 
          where,
          in a nutshell, 
          the \textit{coercivity} inequality is replaced by
          $z\cdot s''(\rho)A(\rho )z \geq \gamma \sum_i
          \rho_i^{2(m-1)} z_i^2$ for some $m$. 
          In that context,
          an $L^2(0,T;H^1(\Omega))$ estimate for $\rho_i$ might be out
          of reach. Instead, 
          the system is rewritten there by using $\rho_i^{m}$ as
          an $L^2(0,T;H^1(\Omega))$ function. Moreover,
          also 
          entropy
          densities~$s$, which are not bounded, such as $s(u)=u-\log u$,
are allowed in~\cite{jungel2015boundedness}. As a consequence,
a different version of the hypothesis (H2b) is considered
  there. We believe that our ansatz can be extended to that 
scenario applying the ideas from \cite{jungel2015boundedness}. This
will mainly affect 
the entropy inequalities and the proof of Proposition \ref{eps.to.0}.
However, we chose to use our simplified assumptions that already cover a large class of parabolic systems. By this, we try to keep the idea of our proof as fundamental as possible.
	\end{remark}
}

\section{Space-time Galerkin method}\label{sec:STGM}
Let the time $T\in(0,\infty)$ be fixed. 
We denote by $Q_T=(0,T)\times\Omega$ the space-time cylinder for a
domain $\Omega\subset\mathbb R^d$, $d\ge 1$.
{We derive our method in four steps.}\medskip

{\bf Step 1 (space-time variational formulation).}
The first step
is to perform integration by parts
in the time variable in~\eqref{eq.vweak}, and to use the
embedding
\begin{equation}\label{eq:C0L2H1embedding}
C([0,T];L^2(\Omega)^N)\subset
L^2(0,T;H^1(\Omega)^N)\cap H^1(0,T;(H^1(\Omega)')^N),
\end{equation}
which can
be proved exactly as in \cite[Chapter 5.9, Theorem 3]{Evans}.
{Then, we define the following variational formulation of
  \eqref{eq.classical}.}

\begin{definition}[{space-time variational formulation/weak
    solution to~\eqref{eq.classical}}]\label{def:weak}
  Find $\rho\in L^2(0,T;H^1(\Omega)^N)\cap
  H^1(0,T;(H^1(\Omega)')^N)$ such that
  \begin{multline}
\label{eq.vweak1}
    \int_{\Omega}\phi(T) \cdot\rho(T)dx
    -\int_0^T\int_{\Omega}\partial_t\phi \cdot\rho dxdt 
\\
+ \sum_{i,j=1}^N\int_0^T\int_{\Omega}\nabla\phi_i \cdot A_{ij}(\rho)\nabla \rho_j dx dt
= \int_0^T\int_{\Omega}\phi \cdot f(\rho) dx dt {+\int_{\Omega}\phi(0) \cdot \rho_0 dx}
\end{multline}
for all $\phi\in H^1(Q_T)^N$.
\end{definition}
  
The following lemma, which will be proven in section
\ref{sect:existence_weak} below (see Remark~\ref{rem:lemma1}),
{establishes that the variational formulation in
  Definition~\ref{def:weak} is actually equivalent to the
  one in~\eqref{eq.vweak}.}

\begin{lemma}\label{lem1}
Let $T>0$. A function $\rho\in L^2(0,T;H^1(\Omega)^N)\cap
H^1(0,T;(H^1(\Omega)')^N)$ {satisfies~\eqref{eq.vweak1} for all $\phi\in H^1(Q_T)^N$ if and only if
it satisfies~\eqref{eq.vweak}} for all $\phi\in L^2(0,T;H^1(\Omega)^N)$, with
$\rho(0)=\rho_0$.
\end{lemma}

{\bf Step 2 (transformation to entropy variables).}
{In the second step, we express
  $\rho$ in 
    formulation~\eqref{eq.vweak1} as a function of the entropy
    variable $w$, namely $\rho=u(w)$. The resulting space-time
    variational formulation is the following:}
{Find $w\in H^1(Q_T)^N$ such that}
\begin{multline}
\label{eq.weak}
\overbrace{
    \int_{\Omega}\phi(T) \cdot {u(w)} (T)dx
    -\int_0^T\int_{\Omega}\partial_t\phi \cdot {u(w)} dxdt }^{{a({u(w)},\phi):=}}
\\
+ \sum_{i,j=1}^N\int_0^T\int_{\Omega}\nabla\phi_i \cdot A_{ij}({u(w)})\nabla ({u(w)})_j dx dt
= \int_0^T\int_{\Omega}\phi \cdot f({u(w)}) dx dt {+\int_{\Omega}\phi(0) \cdot \rho_0 dx}
\end{multline}
for all $\phi\in H^1(Q_T)^N$.
Here, we use the notation $\phi(t):= \mathrm{tr}(\phi)(t,\cdot)$, where $\mathrm{tr}$ denotes the trace operator $\mathrm{tr}:H^1(Q_T)^N\to L^2(\{0,T\}\times\Omega)^N$.\medskip

{\bf Step 3 (regularization).}
Then, we introduce the following
  regularized problem {with regularization parameter
    $\varepsilon$:}
Find ${w^\varepsilon}\in H^1(Q_T)^N$ such that
\begin{multline}
\label{eq.stable}
\varepsilon(\phi,{w^\varepsilon})_{{\WeightedSpace}(Q_T)^N}+a({u(w^\varepsilon)},\phi) 
+ \sum_{i,j=1}^N\int_0^T\int_{\Omega}\nabla\phi_i \cdot A_{ij}({u(w^\varepsilon)})\nabla ({u(w^\varepsilon)})_j dx dt\\
= \int_0^T\int_{\Omega}\phi \cdot f({u(w^\varepsilon)}) dx dt {+\int_{\Omega}\phi(0) \cdot \rho_0 dx}
\end{multline}
for all $\phi\in H^1(Q_T)^N$. 
{The regularization term is given by the scaled
  $H^1(Q_T)^N$ inner product defined as
	\[(f,g)_{H_\varepsilon^1(Q_T)^N}:=\sum_{i=1}^N\left(\int_{Q_T}f_ig_i dxdt + \int_{Q_T} \nabla f_i\cdot \nabla g_i dxdt + \varepsilon\int_{Q_T}\partial_t f_i\partial_t g_i dxdt\right),\] 
    for $f,g\in H^1(Q_T)^N$.  
    Its associated norm will be denoted by $\|\cdot\|_{\WeightedSpace(Q_T)}$.
}
\medskip

{\bf Step 4 (Galerkin discretization).}
Finally, we discretize
equation~\eqref{eq.stable}.
Let $\{\bm V_h\}_{h>0}$ be a family of finite dimensional spaces,
parametrized by $h>0$, such that, for every $h$, $\bm V_h\subset {H^1(Q_T)^N}$.
We make the following approximability
assumption on the family of spaces $\{\bm V_h\}_{h>0}$.
\begin{enumerate}
\item[(H4)]  For all $v\in H^1(Q_T)^N$,
\begin{equation*}
\lim_{h\to 0} \inf_{v_h\in \bm V_h}\|v-v_h\|_{H^1(Q_T)^N}=0.
\end{equation*}
\end{enumerate}

Therefore, we consider the following space-time Galerkin scheme in the entropy
variable unknown:
Find $w_h^\varepsilon\in\bm V_h$ such that
\begin{multline}
\label{eq.galerkin}
\epsilon(\phi,w_h^\varepsilon)_{{\WeightedSpace}(Q_T)^N}+a({u(w_h^\varepsilon)},\phi) 
+\sum_{i,j=1}^N\int_0^T\int_{\Omega}\nabla\phi_i \cdot A_{ij}({u(w_h^\varepsilon)})\nabla ({u(w_h^\varepsilon)})_j dx dt\\
= \int_0^T\int_{\Omega}\phi \cdot f({u(w_h^\varepsilon)}) dx dt {+\int_{\Omega}\phi(0) \cdot \rho_0 dx}
\end{multline}
for all $\phi\in \bm V_h$.

The first term in \eqref{eq.galerkin} can be interpreted as a stabilization term for the Galerkin scheme, with parameter $\varepsilon>0$. 
This is used to obtain a control of the entropy variable. 
Note that, 
due to the nonlinearity of $u$, we expect that
$u(w_h^\varepsilon)\notin \bm V_h$.\medskip

The following two propositions {constitute the main result of this
  paper.
Proposition~\ref{prop.nonlin} establishes that the Galerkin
  problem~\eqref{eq.galerkin} admits solutions $w_h^\varepsilon$ and these solutions
  satisfy an entropy estimate. Then, this is exploited in
  Proposition~\ref{prop2} to obtain existence of weak solutions $\rho$
  to the continuous
  problem~\eqref{eq.classical}, together with related entropy estimates.
Propositions~\ref{prop.nonlin} and~\ref{prop2} will be proven in
section~\ref{sect:existence} and section~\ref{sect:existence_weak},
respectively.
Here and in the following, $|\Omega|$ denotes the volume of $\Omega$,
and
$\gamma$ and $C_f$ are as in Assumption (H2).}

\begin{proposition}[Existence of discrete solutions]\label{prop.nonlin}
There exists a solution $w_h^\varepsilon \in \bm V_h$ of method \eqref{eq.galerkin}. Moreover, every solution  $w_h^\varepsilon \in \bm V_h$ of~\eqref{eq.galerkin}, for $\varepsilon,h>0$, satisfies the entropy estimate 
\begin{equation}\label{w.entropy}
\epsilon\|w_h^\varepsilon
    \|_{{\WeightedSpace}(Q_T)^N}^2+\int_{\Omega}s({u(w_h^\varepsilon (T))})  
    dx+ \gamma\int_{Q_T}|\nabla {u(w_h^\varepsilon)} |^2dxdt\leq \int_{\Omega}s(\rho_0)dx+ C_f|\Omega|T.
  \end{equation}
\end{proposition}

\begin{proposition}[Convergence]\label{prop2}
  Let $w_{h}^\epsilon\in \bm V_h$ be a solution of \eqref{eq.galerkin} for $\epsilon, h>0$. 
Then there exist a weak solution 
\[\rho\in L^2(0,T;H^1(\Omega)^N)\cap H^1(0,T;(H^1(\Omega)')^N)\cap L^\infty((0,T)\times\Omega)^N\]
of \eqref{eq.classical} and sequences $h_i,\epsilon_i\to0$, as $i\to\infty$, such that 
\[u(w^{\epsilon_i}_{h_i})\to\rho\qquad\mbox{in }L^r(Q_T)^N\mbox{, as }i\to \infty\]
for all $r\in[1,\infty)$. Moreover, $\rho$ satisfies the entropy estimate
\begin{equation}\label{eq.entropy.weak}	\int_{\Omega}s(\rho(\tau))dx+\gamma \int_0^\tau\int_{\Omega}\left|\nabla\rho\right|^2 dxdt
\leq  \int_{\Omega}s(\rho_0) dx+C_f|\Omega|\tau
\end{equation}
for all $\tau\in(0,T]$. 
\end{proposition}

{\begin{remark}[Non-closed systems]
{The general setting given by \eqref{eq.classical} describes
  \textit{closed} systems, i.e., without influx or outflux at $\partial\Omega$. These systems are of particular interest, as they obey
  the second law of thermodynamics proving the decay of the
  entropy. However, in some cases (e.g. the lung model presented in
  section \ref{numerical_tests}), one is interested in a non-closed
  subsystem involving, for instance, inhomogeneous Dirichlet boundary
  conditions on some part $\Gamma_D\subset \partial\Omega$ of the
  boundary. If $\rho=g$ is prescribed on $\Gamma_D$, for a given $g\in
  H^1(\Omega)^N$ taking values in $\mathcal{D}$,
we assume that the approximation spaces $\bm V_h$ are affine subspaces of
$s'(g)+ H^1_D(Q_T)^N$, where $H^1_D(Q_T)^N$ is the closure in
$H^1(Q_T)^N$ of the space of $C^{\infty}(Q_T)^N$ functions vanishing at
$\Gamma_D\times(0,T)$. 
Unfortunately, we can neither guarantee the discrete entropy estimate \eqref{w.entropy} nor its continuous version \eqref{eq.entropy.weak}. Instead, we have to work with the relative entropy
			\[
				s^\infty(\rho|g):= s(\rho)-s(g)-s'(g)\cdot(\rho-g),
			\]
                        which is still convex in the variable~$\rho$.
                        We conjecture that,
under the assumption $s'(g)\in L^\infty(\Omega)$,
                        one can prove corresponding versions of Proposition \ref{prop.nonlin} and Proposition \ref{prop2} by employing an estimate of the relative entropy of the form
			\begin{equation*}
				\int_{\Omega}s^\infty(\rho(\tau)|g)dx+\frac\gamma2 \int_0^\tau\int_{\Omega}\left|\nabla\rho\right|^2 dxdt
				\leq  \int_{\Omega}s^\infty(\rho_0|g) dx+(C_f|\Omega|+C_{A,g,f})\tau
			\end{equation*}
for all $\tau\in(0,T)$ and some $C_{A,g,f}$ only depending on $\gamma$,
$\|g\|_{H^1(\Omega)}$ and $\|s'(g)\|_{L^\infty(\Omega)}$, as well as on
$\|A\|_{L^\infty(\mathcal{D})}$ and $\|f\|_{L^\infty(\mathcal D)}$.
This estimate can be derived from~\eqref{eq.vweak}, with test function $\phi=s'(\rho)-s'(g)$. 
However, the proof of convergence of the discrete scheme
\eqref{eq.galerkin}
remains an open problem and is under ongoing investigation.}
\end{remark}	
}

\subsection{Existence of a solution of the numerical scheme}\label{sect:existence}
\begin{proof}[Proof of Proposition \ref{prop.nonlin}]
The idea is to use {the} Leray-Schauder fixed-point theorem for the
    mapping $\Phi:\bm V_h\to \bm V_h$, $v\mapsto {w_h^\varepsilon}$, where ${w_h^\varepsilon}$
    denotes the unique solution of \eqref{eq.galerkin} 
    {with all occurrences of $u(w_h^\varepsilon)$ replaced by $u(v)$.}
    Since $A,f,u$ are continuous, so is
    $\Phi$. Since $\bm V_h$ has finite dimension, $\Phi$ is
      also compact. Then by the Leray-Schauder fixed-point
    theorem, we obtain that $\Phi$ admits a fixed-point if we can
    show that the set
\[\{w\in \bm V_h: w=\sigma\Phi(w),\ \sigma\in[0,1]\}\]
is bounded.

Let ${w_h^\varepsilon}=\sigma \Phi({w_h^\varepsilon})$ for $\sigma\in(0,1]$ and choose $\phi:={w_h^\varepsilon}$. Then \eqref{eq.galerkin} entails
\begin{multline*}
\frac{\epsilon}{\sigma}\| {w_h^\varepsilon}\|_{{\WeightedSpace}(Q_T)^N}^2+
\int_{\Omega}{w_h^\varepsilon} (T) \cdot {u(w_h^\varepsilon (T))} dx
-\int_0^T\int_{\Omega}\partial_t {w_h^\varepsilon} \cdot {u(w_h^\varepsilon)} dxdt 
\\
+ \sum_{i,j=1}^N\int_0^T\int_{\Omega}\nabla ({w_h^\varepsilon})_i \cdot A_{ij}({u(w_h^\varepsilon)})\nabla ({u(w_h^\varepsilon)})_j dx dt
= \int_0^T\int_{\Omega}{w_h^\varepsilon} \cdot f({u(w_h^\varepsilon)}) dx dt+\int_{\Omega}{w_h^\varepsilon} (0) \cdot \rho_0 dx.
\end{multline*}
Using that 
$\partial_t(s(u({w_h^\varepsilon})))=s'(u({w_h^\varepsilon}))\cdot \partial_t(u({w_h^\varepsilon}))={w_h^\varepsilon} \cdot \partial_t(u({w_h^\varepsilon}))$,  we have
\begin{align*}
\partial_t {w_h^\varepsilon} \cdot {u(w_h^\varepsilon)}&=
                                                               \partial_t({w_h^\varepsilon}\cdot u({w_h^\varepsilon}))-{w_h^\varepsilon}\cdot \partial_t(u({w_h^\varepsilon}))=\partial_t ({w_h^\varepsilon}\cdot u({w_h^\varepsilon})-s(u({w_h^\varepsilon}))).
\end{align*}
Thus, by the fundamental theorem of calculus,
\begin{align*}
&\int_{\Omega}{w_h^\varepsilon} (T) \cdot
                 {u(w_h^\varepsilon(T))} dx
                 -\int_{\Omega}{w_h^\varepsilon} (0) \cdot \rho_0 dx
-\int_0^T\int_{\Omega}\partial_t {w_h^\varepsilon} \cdot {u(w_h^\varepsilon)} dxdt
\\&\qquad=
-\int_{\Omega}\big(s({u(w_h^\varepsilon(0))})+{w_h^\varepsilon} (0)\cdot (\rho_0-{u(w_h^\varepsilon(0))})\big)dx
+\int_{\Omega}s({u(w_h^\varepsilon(T))}) dx.
\end{align*}
Note that, by definition, $s'(u({w_h^\varepsilon}))={w_h^\varepsilon}$.
The convexity of $s$ then implies that 
\[s({u(w_h^\varepsilon(0))})+{w_h^\varepsilon} (0)\cdot (\rho_0-{u(w_h^\varepsilon(0))})=s({u(w_h^\varepsilon(0))})+s'({u(w_h^\varepsilon(0))})\cdot (\rho_0-{u(w_h^\varepsilon(0))})\leq s(\rho_0)\]
and hence,
    \begin{align*}
\int_{\Omega}{w_h^\varepsilon} (T) \cdot {u(w_h^\varepsilon(T))} dx-\int_{\Omega}{w_h^\varepsilon} (0) \cdot \rho_0 dx
-\int_0^T\int_{\Omega}\partial_t {w_h^\varepsilon} \cdot {u(w_h^\varepsilon)} dxdt\\
\geq
\int_{\Omega}s({u(w_h^\varepsilon(T))}) dx
-	\int_{\Omega}s(\rho_0)dx.
\end{align*}
The next step is to use (H2a) in combination with ${w_h^\varepsilon}=s'({u(w_h^\varepsilon)})$, which yields that 
\begin{align*}
    \sum_{i,j=1}^N\nabla ({w_h^\varepsilon})_i\cdot A_{ij}({u(w_h^\varepsilon)})\nabla ({u(w_h^\varepsilon)})_j
      =\sum_{i,j=1}^N \nabla (s'({u(w_h^\varepsilon)}))_i 
      \cdot A_{ij}({u(w_h^\varepsilon)})\nabla ({u(w_h^\varepsilon)})_j\\
      = \sum_{i,j,k=1}^N\nabla ({u(w_h^\varepsilon)})_k\cdot (s''({u(w_h^\varepsilon)}))_{ki}A_{ij}({u(w_h^\varepsilon)})\nabla ({u(w_h^\varepsilon)})_j\geq \gamma |\nabla {u(w_h^\varepsilon)}|^2,
\end{align*}
where $|\nabla {u(w_h^\varepsilon)}|^2:=\sum_{\ell=1}^d |\frac{\partial}{\partial x_\ell}{u(w_h^\varepsilon)}|^2$.
Moreover, due to (H2b)
and ${w_h^\varepsilon}=s'({u(w_h^\varepsilon)})$, we have
\begin{align*}
{w_h^\varepsilon}\cdot f({u(w_h^\varepsilon)})=
s'({u(w_h^\varepsilon)})\cdot f({u(w_h^\varepsilon)})\leq C_f .
\end{align*}
Therefore, we can conclude the entropy estimate
\begin{equation*}
\frac{\epsilon}\sigma \|{w_h^\varepsilon}\|_{{\WeightedSpace}(Q_T)^N}^2+\int_{\Omega}s({u(w_h^\varepsilon(T))})dx+\gamma \int_{Q_T}|\nabla {u(w_h^\varepsilon)}|^2dxdt\leq \int_{\Omega}s(\rho_0)dx+ C_f|\Omega|T. 
\end{equation*}
Hence, $\|{w_h^\varepsilon}\|_{{\WeightedSpace}(Q_T)^N}^2$ is uniformly bounded, because $\sigma\leq 1$. 
Thus, the Leray-Schauder theorem is applicable and yields that $\Phi$ has a fixed point, and therefore the scheme \eqref{eq.galerkin} admits a solution.
Using these calculations for $\sigma=1$, it follows that every
    solution has to satisfy the 
    entropy inequality~\eqref{w.entropy}.
  \end{proof}
  
\subsection{Convergence of the numerical scheme as \texorpdfstring{$h\to0$}{h->0}}\label{sec:limh}
We will show that, for a fixed $\epsilon>0$, the numerical scheme~\eqref{eq.galerkin} converges as $h\rightarrow 0$.

\begin{proposition}[Convergence of the scheme for fixed $\epsilon>0$]\label{prop3}
There exists ${w^\varepsilon}\in H^1(Q_T)^N$ with ${\rho^\varepsilon} := u({w^\varepsilon})\in
L^2(0,T,H^1(\Omega)^N)$,
and a sequence $h_\ell\to0$  such that
\[
{\rho_{h_\ell}^\varepsilon}:= u({w_{h_\ell}^\varepsilon})\to {\rho^\varepsilon}\mbox{ strongly in } L^r(Q_T) \mbox{ for all }r\in[1,\infty).\]
Moreover, ${w^\varepsilon}$ {solves} 
\eqref{eq.stable} and {satisfies} 
the entropy estimate
\begin{equation}\label{eq.entropy-cont}
\epsilon\|{w^\varepsilon} \|_{{\WeightedSpace}(Q_T)^N}^2+\int_{\Omega}s({\rho^\varepsilon} (T))  
dx+ \gamma\int_{Q_T}|\nabla {\rho^\varepsilon}|^2dxdt\leq \int_{\Omega}s(\rho_0)dx+ C_f|\Omega|T.
\end{equation}
\end{proposition}
\begin{proof}
The first part of the assertion follows from the fact that
    ${w_{h}^\varepsilon}$ is uniformly bounded in ${\WeightedSpace}(Q_T)^N$, which yields
    that there exists ${w^\varepsilon}\in H^1(Q_T)^N$ and  subsequence
    ${h_\ell}\to 0$ such that ${w_{h_\ell}^\varepsilon}\rightharpoonup {w^\varepsilon}$ in
    ${\WeightedSpace}(Q_T)^N$, due to the Banach-Alaoglu theorem, and
    ${w_{h_\ell}^\varepsilon}\to {w^\varepsilon}$ in $L^2(Q_T)^N$, due to Rellich's theorem. In particular, we can choose this subsequence in such a way that ${w_{h_\ell}^\varepsilon}$ converges a.e.\ to ${w^\varepsilon}$.
As $u$ is bounded (see Assumption~(H2)), the dominated
convergence theorem entails the strong
convergence of ${\rho_{h_\ell}^\varepsilon}\equiv u({w_{h_\ell}^\varepsilon})\to u({w^\varepsilon})=:{\rho^\varepsilon}$
in $L^r(Q_T)^N$ for all $r\in[1,\infty)$.
Combining this with the entropy estimate \eqref{w.entropy}, there exists another subsequence (which we do not relabel) such that ${\rho_{h_\ell}^\varepsilon}\rightharpoonup {\rho^\varepsilon}$ weakly in $L^2(0,T,H^1(\Omega)^N)$.

Finally, owing to assumption (H4), for every $\phi\in H^1(Q_T)^N$, there exists $\phi_{h_\ell}\in \bm V_{h_\ell}$ such that $\phi_{h_\ell}\to \phi$ in ${\WeightedSpace}(Q_T)^N$. 
Using $\phi_{h_\ell}$ as a test function in \eqref{eq.galerkin}, we obtain \eqref{eq.stable} in the limit $h_\ell\to0$, as each integral in \eqref{eq.galerkin} converges separately.
The entropy inequality~\eqref{eq.entropy-cont}
is a consequence of Fatou's lemma and the weak lower semicontinuity of the norm.
\end{proof}

The following corollary will be used in the analysis of the limit
for $\epsilon\to 0$ (see proof of Proposition~\ref{eps.to.0} below).

\begin{corollary}\label{cor1}
  Let $\tau,\delta\geq0$ be such that $\tau+\delta\leq T$.
It holds true that
\begin{multline}\label{eq.entropy.reg42b}	\epsilon \| {w^\varepsilon}\|_{{\WeightedSpace}(Q_\tau)^N}^2+\frac1\delta\int_{\tau}^{\tau+\delta}\!\int_{\Omega}s({\rho^\varepsilon})dxdt
+ \gamma\left(1+{\frac{\sqrt\varepsilon}{\delta}}\right)\int_0^\tau\int_{\Omega}|\nabla {\rho^\varepsilon}|^2 dxdt
\\\leq  	\left(1+{\frac{\sqrt\varepsilon}{\delta}}\right)\int_{\Omega}s(\rho_0) dx
+C_f|\Omega|\left(\tau+{\frac{\delta}{2}+\frac{\sqrt\varepsilon}{\delta}T}\right), 
  \end{multline}
where {$\rho^\varepsilon := u(w^\varepsilon)$ and} $Q_\tau:=(0,\tau)\times\Omega$.
\end{corollary}
\begin{proof}
Set 
\begin{equation*}
\psi(t):=\begin{cases}
      1 & \text{if}\ t<\tau,\\
      1-\frac{t-\tau}\delta &\text{if}\ \tau\leq t \leq \tau+\delta,
\\
0 &\mbox{otherwise}.
\end{cases}
\end{equation*}
Thus, ${w^\varepsilon}\psi\in H^1(Q_T)^N$. 
 Similarly as in the proof of Proposition \ref{prop.nonlin}, we use {$\rho^\varepsilon := u(w^\varepsilon)$} and 
\begin{align*}
\partial_t (\psi {w^\varepsilon}) \cdot {\rho^\varepsilon}= \partial_t(\psi {w^\varepsilon}\cdot {\rho^\varepsilon})-\psi {w^\varepsilon}\cdot \partial_t {\rho^\varepsilon}=\partial_t (\psi {w^\varepsilon}\cdot {\rho^\varepsilon}-\psi s({\rho^\varepsilon}))+\partial_t\psi s({\rho^\varepsilon})
\end{align*}
and, since $\psi(T)=0$ and $\psi(0)=1$,
\begin{align*}
    \begin{split}
      \int_{Q_T}\!\!\!\partial_t ({w^\varepsilon}\psi)  \cdot {\rho^\varepsilon} 
      dxdt+\int_{\Omega}\!\!{w^\varepsilon} (0)\cdot \rho_0 dx
=\int_{Q_T}\!\!\!\partial_t \psi s({\rho^\varepsilon}) dxdt
+
\int_{\Omega}\!\!\big(s({\rho^\varepsilon} (0))+{w^\varepsilon} (0)\cdot (\rho_0-{\rho^\varepsilon} (0))\big)dx.
\end{split}
\end{align*}
Thus, using the definition of $\psi$, and treating the last term of the previous equation as in the proof of Proposition~\ref{prop.nonlin}, we get
\begin{equation*}
\int_{Q_T}\partial_t ({w^\varepsilon}\psi)  \cdot {\rho^\varepsilon} dxdt+\int_{\Omega}\psi(0) {w^\varepsilon} (0)\cdot \rho_0 dx+\frac{1}{\delta}\int_\tau^{\tau+\delta}\!\int_{\Omega} s({\rho^\varepsilon}) dxdt\leq \int_{\Omega}s(\rho_0) dx.
  \end{equation*}
From \eqref{eq.stable} tested with $\phi=\psi{w^\varepsilon}$ and the previous
inequality, we get
\begin{align*}
\begin{split}
\epsilon(\psi {w^\varepsilon}, {w^\varepsilon})_{{\WeightedSpace}(Q_T)^N}
+\sum_{i,j=1}^N\int_0^T\int_{\Omega}\nabla (\psi {w^\varepsilon})_i \cdot
A_{ij}({\rho^\varepsilon})\nabla ({\rho^\varepsilon})_j dx dt
+\frac{1}{\delta}\int_\tau^{\tau+\delta}\!\int_{\Omega} s({\rho^\varepsilon})
dxdt\\
\le
\int_{\Omega}s(\rho_0) dx
+\int_0^T\int_{\Omega}\psi {w^\varepsilon} \cdot f({\rho^\varepsilon}) dx dt
\end{split}
\end{align*}
which, due to the properties of $\psi$ and the assumption (H2), entails      
  \begin{multline*}	\epsilon(\psi {w^\varepsilon}, {w^\varepsilon})_{{\WeightedSpace}(Q_T)^N}
    +\frac1\delta\int_{\tau}^{\tau+\delta}\!\int_{\Omega}s({\rho^\varepsilon})dxdt
    + \gamma\int_0^\tau\int_{\Omega}|\nabla {\rho^\varepsilon}|^2 dxdt
  \leq  \int_{\Omega}s(\rho_0) dx
   +C_f|\Omega|(\tau+\delta/2).
\end{multline*}
Finally, we can estimate the first term as
\begin{align*}
\epsilon(\psi {w^\varepsilon}, {w^\varepsilon})_{{\WeightedSpace}(Q_T)^N}&=
\epsilon\int_0^T\int_{\Omega}\left(
\psi |{w^\varepsilon}|^2+\psi |\nabla {w^\varepsilon}|^2+ {\varepsilon}\partial_t(\psi {w^\varepsilon})\cdot\partial_t {w^\varepsilon}
\right) dxdt\\\notag
&=\epsilon\int_0^T\int_{\Omega}\left(\psi |{w^\varepsilon}|^2+\psi |\nabla {w^\varepsilon}|^2+{\varepsilon}\psi
|\partial_t {w^\varepsilon}|^2+{\varepsilon}\partial_t\psi {w^\varepsilon}\cdot \partial_t {w^\varepsilon}
\right) dxdt\\ 
&\ge \epsilon \| {w^\varepsilon}\|_{{\WeightedSpace}(Q_\tau)^N}^2
-{\frac{\epsilon^2}{\delta}}\int_\tau^{\tau+\delta}\int_{\Omega} {w^\varepsilon}\cdot\partial_t {w^\varepsilon} dxdt.
\end{align*}
Using the Cauchy-Schwarz inequality and the definition of
      the ${\WeightedSpace}$ norm
    yields
\[
{\frac{\epsilon^2}{\delta}}\int_\tau^{\tau+\delta}\int_{\Omega} {w^\varepsilon}\cdot\partial_t {w^\varepsilon} dxdt
\leq 
{
	\frac{\epsilon^2}{\delta} \|{w^\varepsilon}\|_{L^2((\tau,\tau+\delta)\times\Omega)^N}\|\partial_t{w^\varepsilon}\|_{L^2((\tau,\tau+\delta)\times\Omega)^N}
\leq \frac{\epsilon^{\frac32}}{\delta} \|{w^\varepsilon}\|^2_{\WeightedSpace((\tau,\tau+\delta)\times\Omega)^N}
},
\]
and therefore
\begin{multline*}
\epsilon \| {w^\varepsilon}\|_{{\WeightedSpace}(Q_\tau)^N}^2+\frac1\delta\int_{\tau}^{\tau+\delta}\!\int_{\Omega}s({\rho^\varepsilon})dxdt
+\gamma\int_0^\tau\int_{\Omega}|\nabla {\rho^\varepsilon}|^2 dxdt
\\  \leq {\frac{\epsilon^\frac32}{\delta}} \|{w^\varepsilon}\|_{{H^1_{\varepsilon}}((\tau,\tau+\delta)\times\Omega)^N}^2+  \int_{\Omega}s(\rho_0) dx
+C_f|\Omega|(\tau+\delta/2).
\end{multline*}
Note that we cannot estimate the first term on the right-hand side by
the first term on the left-hand side, because the domain of the
norms are disjoint.
Fortunately, we have the entropy estimate \eqref{eq.entropy-cont}, which we add {$\sqrt\epsilon/\delta$} times to this inequality to get
\begin{multline*}
\epsilon\| {w^\varepsilon}\|_{{\WeightedSpace}(Q_\tau)^N}^2+\frac1\delta\int_{\tau}^{\tau+\delta}\!\int_{\Omega}s({\rho^\varepsilon})dxdt
+{\frac{\sqrt\varepsilon}{\delta}}\int_{\Omega}s({\rho^\varepsilon} (T))  
dx
+\gamma(1+{\frac{\sqrt{\varepsilon}}{\delta}})\int_0^\tau\int_{\Omega}|\nabla {\rho^\varepsilon}|^2 dxdt
\\  \leq 	(1+{\frac{\sqrt\varepsilon}{\delta}})\int_{\Omega}s(\rho_0) dx
+C_f|\Omega|\left(\tau+{\frac{\delta}{2}+\frac{\sqrt\varepsilon}{\delta}T}\right). 
\end{multline*}
which, since $s({\rho^\varepsilon} (T))\ge 0$, implies the assertion.
\end{proof}

\subsection{Limit of \texorpdfstring{$\epsilon\to0$}{eps->0}}\label{sec:limeps}
We consider the limiting problem
\begin{multline}
\label{eq.limit}
-\int_{\Omega}\phi(0) \cdot \rho_0 dx
-\int_0^T\int_{\Omega}\partial_t\phi \cdot \rho dxdt 
+ \sum_{i,j=1}^N\int_0^T\int_{\Omega}\nabla\phi_i \cdot A_{ij}(\rho)\nabla \rho_j dx dt
\\= \int_0^T\int_{\Omega}\phi \cdot  f(\rho) dx dt
\end{multline}
for all $\phi\in (H^1(Q_T))^N$ with $\phi(T)=0$. As above,
we use the notation $\phi(t):= \mathrm{tr}(\phi)(t,\cdot)$, where $\mathrm{tr}$ denotes the trace operator $\mathrm{tr}:H^1(Q_T)^N\to L^2(\{0,T\}\times\Omega)^N$.
\begin{proposition}\label{eps.to.0}
  Let $\tau,\delta\geq0$ such that $\tau+\delta\leq T$.
{Set $\rho^\varepsilon:=u(w^\varepsilon)$.} Then there exist $\rho\in L^2(0,T;H^1(\Omega)^N)$ with $\rho(t,x)\in\overline{\mathcal{D}}$ for a.e.~$(t,x)\in Q_T$ being a solution of \eqref{eq.limit} and a subsequence $\epsilon_j\to0$ such that
\[\rho^{\epsilon_j}\to \rho\qquad\mbox{in every }L^r(Q_T)^N, r\in[1,\infty),\mbox{ as }\epsilon_j\to0. \]
Moreover, $\rho$ satisfies the entropy inequality

\begin{equation}\label{entropy42}
\frac1\delta\int_{\tau}^{\tau+\delta}\!\int_{\Omega}s(\rho)dxdt
+ \gamma\int_0^\tau\int_{\Omega}|\nabla \rho|^2 dxdt
\leq {\int_{\Omega}s(\rho_0) dx} + C_f|\Omega|
{(\tau+\delta/2)}.
\end{equation}
\end{proposition}

In the proof of Proposition~\ref{eps.to.0}, the key ingredient to prove strong convergence of (at least a subsequence of) $\rho^\epsilon$
will be the idea of compensated compactness, which is a special
technique applying the classical div-curl lemma; see,
e.g.~{\cite[Lemma~7.2]{TartarBook2009}}. 

\begin{lemma}[div-curl lemma]
Let
    $\alpha,\alpha^\ell\in L^2(Q_T)^{1+d}$ and $\beta,\beta^\ell\in L^2(Q_T)^{1+d}$.
Then
\begin{align*}
&\alpha^\ell\rightharpoonup \alpha\quad\mbox{in
}L^2(Q_T)^{1+d}\text{ as } \ell\to+\infty,\
\mbox{and}\quad (\mathrm{div}_{(t,x)}\alpha^{\ell})_{\ell\in\mathbb
N}\ \mbox{is bounded in } L^2(Q_T),
\\
&\beta^\ell\rightharpoonup \beta\quad\mbox{in
}L^2(Q_T)^{1+d} \text{ as } \ell\to+\infty,\
\mbox{and}\quad  (\mathrm{curl}_{(t,x)}\beta^{\ell})_{\ell\in\mathbb N}\mbox{ is bounded in } L^2(Q_T)^{(1+d)\times(1+d)}
\end{align*}
implies that 
\begin{align*}
\alpha^\ell\cdot \beta^\ell \rightharpoonup \alpha\cdot \beta
      \qquad \mbox{in }\mathcal D'(Q_T)\quad
      \text{ as } \ell\to+\infty,
    \end{align*}
where $\mathcal D'(Q_T)$ denotes the dual space of $\mathcal D(Q_T):=C_c^\infty(Q_T)$.
\end{lemma}

\begin{proof}[Proof of Proposition \ref{eps.to.0}]
Let $w^\epsilon, \rho^\epsilon:=u(w^\epsilon)$ denote the solution of~\eqref{eq.stable}
satisfying the entropy inequality~\eqref{eq.entropy-cont}.
For any fixed $i$, $i=1,\ldots,N$, we define
the vector-valued functions with $(1+d)$ components
\begin{equation*}
\alpha^\epsilon=\binom{\rho_i^\epsilon-{\epsilon^2}\partial_t w_i^\epsilon}{J_i^\epsilon-\epsilon\nabla w^\epsilon_i}\quad\mbox{and}\quad \beta^\epsilon:= \binom{\rho^\epsilon_i}{0},\qquad \mbox{where }J_i^\epsilon=-\sum_{j=1}^N A(\rho^\epsilon)_{ij} \nabla\rho_j^\epsilon.
\end{equation*}
Note that, by assumption, $\mathcal{D}$ is bounded and so is
$\rho^\epsilon= u(w^\epsilon)$. Thus, thanks to the entropy estimate~\eqref{eq.entropy-cont},
$\alpha^\epsilon,\beta^\epsilon$ are bounded uniformly in $L^2(Q_T)^{1+d}$
w.r.t.~$\epsilon\in(0,1)$. By the Banach-Alaoglu theorem, there exist
$\alpha,\beta\in
L^2(Q_T)^{1+d}
$ and a subsequence $\epsilon_\ell\to0$ such that
\[\alpha^{\epsilon_{\ell}}\rightharpoonup \alpha,\quad
\beta^{\epsilon_{\ell}}\rightharpoonup \beta\quad\mbox{in
}
L^2(Q_T)^{1+d}
\qquad\mbox{as }\epsilon_{\ell}\to0.\]
Clearly, $\beta$ has the form $(\rho_i,0)$ for some $\rho_{i}\in L^2(Q_T)$. 
Due to the entropy estimate~\eqref{eq.entropy-cont},
{$\sqrt{\epsilon}\|w^\epsilon\|_{\WeightedSpace(Q_T)}$}
is bounded. 
Hence,
$\beta^\epsilon_0-\alpha^\epsilon_0={\epsilon^2}\partial_t
w^\epsilon_i\to0$ in $L^2(Q_T)$ as $\epsilon\to0$, implying that
$\rho_i:=\beta_0=\alpha_0$ and $\alpha\cdot\beta=\rho_i^2$, where in
this context the index 0 denotes the first 
component of the $(1+d)$-dimensional vector.
Moreover, one can easily show that
\[\|\mathrm{curl}_{(t,x)}
\beta^\epsilon\|_{L^2(Q_T)^{(1+d)\times(1+d)}}\leq C \|\nabla
\rho_i^\epsilon\|_{L^2(Q_T)^d}\]
for some $C>0$.
Again by the entropy estimate~\eqref{eq.entropy-cont},
this implies that $\mathrm{curl}_{(t,x)} \beta^\epsilon$ is uniformly
bounded\footnote{{The fact that the $L^2$ norm of $\nabla\rho^\epsilon_i$ is uniformly bounded is ultimately a consequence of our hypothesis (H2a), i.e., the matrix $s''(\rho)A(\rho)$ being coercive. Using instead the original assumptions made by A.~J\"ungel in \cite{jungel2015boundedness}, one can only assure that $\nabla \rho_i^{m}$ is bounded in $L^2$ for some $m$. However, one can circumvent this issue by defining $\beta^\epsilon:=(\rho_i^m,0)^T$ in this case.}} in
$L^2(Q_T)^{(1+d)\times(1+d)}$ w.r.t.~$\epsilon\in(0,1)$. 
In order to apply the div-curl lemma, it remains to prove that the
space-time divergence of $\alpha^\epsilon$ is bounded. 
For this, we require the equation for
$\rho_i^\epsilon$ in the
interior of $Q_T$, i.e., 
from equation~\eqref{eq.stable},
\begin{multline*}
\epsilon\int_{Q_T} \psi w^\epsilon_idxdt
+{\epsilon^2}\int_{Q_T} \partial_t\psi
\partial_t w^\epsilon_idxdt
+\epsilon\int_{Q_T} \nabla\psi\cdot \nabla w^\epsilon_idxdt
-\int_{Q_T}\partial_t\psi \rho^\epsilon_i dxdt 
\\
+\sum_{j=1}^N\int_{Q_T}\nabla\psi \cdot A_{ij}(\rho^\epsilon)\nabla \rho^\epsilon_j dx dt
=\int_{Q_T}\psi  f_i(\rho^\epsilon) dx dt
\end{multline*}
for all $\psi\in H^1_0(Q_T)$. We can rewrite this by using the weak space-time divergence of $\alpha^\epsilon$ as
\begin{align*}
-\int_{Q_T}\nabla_{(t,x)}\psi\cdot \alpha^\epsilon dxdt  &=
\int_{Q_T} \partial_t\psi ({\epsilon^2}\partial_t w^\epsilon_i-\rho^\epsilon_i) dxdt 
\\ & \qquad+
\int_{Q_T} \nabla\psi\cdot \left(\epsilon\nabla w^\epsilon_i
+\sum_{j=1}^N A_{ij}(\rho^\epsilon)\nabla \rho^\epsilon_j\right) dx dt
\\
&=\int_{Q_T}\psi  f_i(\rho^\epsilon) dx dt-\epsilon\int_{Q_T} \psi w^\epsilon_idxdt
\end{align*}
for all  $\psi\in H^1_0(Q_T)$. We observe that the right-hand side
defines a bounded operator on $L^2(Q_T)$ due to the entropy
estimate~\eqref{eq.entropy-cont} and the fact that $f_i$ is
    uniformly bounded as a continuous function defined on a
    compact set (see (H2)).
This yields that $\mathrm{div}_{(t,x)}\alpha^\epsilon$ is uniformly bounded in $L^2(Q_T)$. 
Therefore, we can apply the div-curl lemma and obtain
that
\[(\rho^{\epsilon_\ell}_i-{\epsilon_\ell^2}\partial_tw_i^{\epsilon_{\ell}})\rho^{\epsilon_\ell}_i=\alpha^{\epsilon_\ell}\cdot\beta^{\epsilon_\ell}\rightharpoonup
\alpha\cdot\beta =\rho_i^2\quad \mbox{in }\mathcal{D}'(Q_T)\qquad
\text{as } \epsilon_\ell\to0.\]
Using that
$\rho_i^{\epsilon_{\ell}}\rightharpoonup \rho_{i}$ and ${\epsilon_{\ell}^2} {\partial_t}w_i^{\epsilon_{\ell}}\to0$  in $L^2(Q_T)$, we obtain that 
\[\int_{Q_T}(\rho_{i}^{\epsilon_\ell})^2\phi^2 dxdt\to\int_{Q_T}\rho_i^2 \phi^2 dxdt\qquad\mbox{as }\epsilon_{\ell}\to0\]
for all $\phi\in C_c^\infty(Q_T)$. Hence, $\phi\rho_{i}^{\epsilon_{\ell}}\to \phi\rho_i$ in $L^2(Q_T)$ for all $\phi\in C_c^\infty(Q_T)$. In particular, there exists a subsequence not being relabeled such that $\rho^{\epsilon_{\ell}}_i\to\rho_i$ a.e.\ in $Q_T$. For almost every $(t,x)\in Q_T$, we know that $\rho^{\epsilon_{\ell}}(t,x)\in \mathcal{D}$ and that $\mathcal{D}$ is bounded. Thus, we can apply the dominated convergence theorem, which yields that
\[\rho{_i}^{\epsilon_{\ell}}\to \rho{_i}\qquad\mbox{in every }L^r(Q_T), r\in[1,\infty),\mbox{ as }\epsilon_{\ell}\to0, \]
and that $\rho(t,x)\in \mathcal D$ for almost every $(t,x)\in Q_T$.

Moreover, the entropy inequality~\eqref{eq.entropy-cont}
also states that $\nabla\rho_i^\epsilon$ is bounded in
$L^2(Q_T)^d$ independently of $\epsilon$. Since
$|\rho^\epsilon|= |u(w^\epsilon)|=|(s^\prime)^{-1}(w^\epsilon)|\le
\sup_{v\in\mathcal D}|v|^2$, according to (H2), then, using
again~\eqref{eq.entropy-cont}, we obtain
\begin{align*}
\|\rho_i^\epsilon\|_{L^2(0,T;H^1(\Omega))}^2 
&=
\int_{Q_T}(\rho_i^\epsilon)^2dxdt +
\int_{Q_T}|\nabla\rho_i^\epsilon|^2dxdt
\\
&\leq
 |\Omega|T \|\rho_i^\epsilon\|_{L^\infty(Q_T)}^2    
 + \frac1\gamma\left(\int_{\Omega}s(\rho_0)dx+C_f|\Omega|T\right) \\
&\leq \frac{1}{\gamma}\int_{\Omega}s(\rho_0)dx + \left(\sup_{v\in\mathcal D}|v|^2+\frac{C_f}{\gamma}\right)|\Omega|T,
\end{align*}   
namely, $\rho_i^\epsilon$ is bounded in $L^2(0,T;H^1(\Omega))$
independent on $\epsilon$.
Taking yet another subsequence, which we do not relabel, we
can see that there exists $\rho_i\in L^2(0,T;H^1(\Omega))$
such 
that $\rho^{\epsilon_\ell}_i\rightharpoonup \rho_i$ in $L^2(0,T;H^1(\Omega))$. In particular,
$\nabla\rho^{\epsilon_\ell}_i\rightharpoonup\nabla \rho_i$ in
$L^2(Q_T)^d$.
We already have seen that $\sqrt{\epsilon}{\|w^\epsilon\|_{\WeightedSpace(Q_T)}}$ is bounded, then
$\epsilon
{\|w^\epsilon\|_{\WeightedSpace(Q_T)}}\to0$ {implying $\epsilon
{(w^\epsilon,\phi)_{\WeightedSpace(Q_T)}}\to0$}.

Now, we prove that $\rho$ is solution to the limiting
problem~\eqref{eq.limit}.
Let $\phi\in H^1(Q_T)$ with trace $\phi(T)=0$. Using that $A$ is
bounded, according to (H1), the dominated convergence theorem
yields
\[
\int_{Q_T}|\nabla\phi|^2|A_{ij}(\rho^{\epsilon_\ell})|^2dxdt \to
\int_{Q_T}|\nabla\phi|^2|A_{ij}(\rho)|^2dxdt\qquad
\text{as } \epsilon_\ell\to 0.
\]
In particular, $\nabla\phi A_{ij}(\rho^{\epsilon_\ell})$ converges
strongly in $L^2(Q_T)^d$.
For each $i=1,\ldots,N$, we test the equation for
$\rho_i^\epsilon$ (see~\eqref{eq.stable}) with functions
$\phi\in H^1(Q_T)$ with trace $\phi(T)=0$,
take the limit for $\epsilon=\epsilon_{\ell}\to 0$,
and obtain
\begin{multline*}
-\int_{\Omega}\phi(0)  \rho_i^0 dx
-\int_0^T\int_{\Omega}\partial_t\phi \rho_i dxdt 
+ \sum_{j=1}^N\int_0^T\int_{\Omega}\nabla\phi \cdot A_{ij}(\rho)\nabla \rho_j dx dt
\\= \int_0^T\int_{\Omega}\phi  f_i(\rho) dx dt
\end{multline*}
for all $i=1,\ldots,N$.

Finally, recall that $\rho^\epsilon$
satisfies the
entropy estimate~\eqref{eq.entropy.reg42b} from Corollary \ref{cor1}.
Thus, we obtain the entropy inequality~\eqref{entropy42}
as a direct consequence of the lower weak continuity of the $L^2$ norm and the Fatou lemma.
\end{proof}

\subsection{Existence of a weak solution}\label{sect:existence_weak}
In this section, we prove that problem~\eqref{eq.classical}
possesses a weak solution $\rho$ {in the sense of Definition~\ref{def:weak}}. Moreover, we prove the equivalence stated in
Lemma~\ref{lem1} {between} the weak
formulation~\eqref{eq.vweak1} in Definition~\ref{def:weak} and the
weak formulation~\eqref{eq.vweak}.

\begin{proposition}\label{prop:existence_weak}
Let $\rho$ be given by Proposition \ref{eps.to.0}. Then $\rho\in
H^1(0,T;(H^1(\Omega)')^N)$ and $\rho\in C^0([0,T];L^2(\Omega))$ with $\rho(0)=\rho_0$. Moreover, it satisfies the entropy inequality
\begin{equation}\label{eq.entropy.reg73}	\int_{\Omega}s(\rho(\tau))dx+\gamma \int_0^\tau\int_{\Omega}\left|\nabla\rho\right|^2 dxdt
\leq  \int_{\Omega}s(\rho_0) dx+C_f|\Omega|\tau.
\end{equation}
for almost all $\tau\in(0,T)$.
\end{proposition}
\begin{proof}
Using the equation \eqref{eq.limit}, we obtain that 
\begin{align*}
\left|\int_{Q_T}\partial_t\phi \rho_i dxdt \right|
&\leq \sum_{j=1}^N\int_{Q_T}|\nabla\phi|| A_{ij}(\rho)\nabla
\rho_j |dx dt+ \int_{Q_T}|\phi|  |f_i(\rho)| dx
dt+\int_\Omega|\phi(0)||\rho_{0,i}| dx
\\&\leq C_\rho\|\phi\|_{{L^2(0,T;H^1(\Omega))}},
\end{align*}
since $\rho\in L^\infty(Q_T)\cap L^2(0,T;H^1(\Omega))$.
This implies that, for each $i=1,\ldots,N$, $\rho_i$ has a weak time derivative satisfying
    $\partial_t\rho_i\in L^2(0,T;H^1(\Omega)')$.
    Then the embedding $H^1(0,T;H^1(\Omega)')\cap
    L^2(0,T;H^1(\Omega))
    \subset C^0([0,T];L^2(\Omega))$, entails that every
      $\rho_i$ is continuous in time, and so is $\rho$. We obtain the desired entropy estimate as a limit $\delta\to0$ of \eqref{entropy42}.

It remains to show that $\rho(0)= \rho_0$ in
    $L^2(\Omega)^N$. For this, let $\psi\in H^1(\Omega)^N$
    and, for $\tau\in(0,T)$, define
\begin{equation*}
\phi_\tau(t,\cdot):=\begin{cases}
      \left(1-\frac{t}{\tau}\right)\psi(\cdot)&\text{in}\ \Omega\times[0,\tau],\\
0&\text{in}\ \Omega\times(\tau,+\infty).
\end{cases}
\end{equation*}
We easily see that $\phi_\tau\to0$ in
    $L^2(0,T;H^1(\Omega)^N)$ as $\tau\to0$.
    Then, from equation \eqref{eq.limit} tested with
      $\phi_\tau$, we get
    that, for all $\psi\in H^1(\Omega)^N$,
\begin{align*}
\int_{\Omega}\left(\frac{1}{\tau}\int_0^\tau \rho
      dt-\rho_0\right)\psi dx\to 0\qquad \text{as $\tau\to0$.}
\end{align*}
 Finally, the continuity implies that $\lim_{\tau\to0}\frac{1}{\tau}\int_0^\tau \rho dt=\rho(0)$, which entails $\rho(0)=\rho_0$.
\end{proof}

\begin{remark}\label{rem:lemma1}
Using the last part of the proof of
Proposition~\ref{prop:existence_weak}, we can easily show that any
solution $\rho$ {of \eqref{eq.vweak1}} 
satisfies
$\rho(0)=\rho_0$. Therefore, the proof of Lemma~\ref{lem1} is a straightforward
application of the integration by parts formula and of the embedding~\eqref{eq:C0L2H1embedding}.
\end{remark}

\begin{corollary}
Let $\rho$ be given by Proposition \ref{eps.to.0}. Then $\rho$
    is a solution of \eqref{eq.vweak}. 
\end{corollary}
\begin{proof}
Thanks to Proposition~\ref{prop:existence_weak}, 
we know that $\rho$ possesses enough regularity such that we can integrate in 
\eqref{eq.limit} w.r.t.~$t$, which yields \eqref{eq.vweak} for all $\phi\in H^1(Q_T)^N$ with $\phi(T)=0$. Using a density argument yields the assertion.
  \end{proof}

The proof of {Proposition}~\ref{prop2} is now straightforward.
  
\begin{proof}[Proof of {Proposition} \ref{prop2}]
We only have to collect the previous results to obtain the proposition using a diagonal sequence argument.
\end{proof}

\section{Applications and numerical tests}\label{sec:numerics}

In this section, we apply the general setting of
section~\ref{sec:generalsetting} and numerically test
the space-time Galerkin method of section~\ref{sec:STGM} by
considering four problems:
the (linear) heat equation (section~\ref{sec:heat}),
the porous medium equation (section~\ref{sec:porousmedium}),
the Fisher-KPP equation (section~\ref{sec:FisherKPP}),
and 
the Maxwell-Stefan system {in the case of $N=2$ species}
(section~\ref{sec.ThreeGas}).
For the Maxwell-Stefan system, the discussion of the general setting {and of an alternative space-time
  Galerkin method for the case of $N>2$}
is postponed to section~\ref{sec.MaxwellStefan}.
We remark that we apply this nonlinear setting to the linear heat
equation for validation purposes and, in particular, in order to
stress its
unconditional stability on a simple test problem.

In all cases,
we consider the entropy density
$s:\mathcal D\to [0,+\infty)$
defined by
\begin{equation}\label{eq:definition_of_s}
s(\rho)=\sum_{j=1}^N\rho_j\log \rho_j+\left(1-\sum_{j=1}^N\rho_j\right)\log\left(1-\sum_{j=1}^N\rho_j\right)+\log(N+1),
\end{equation}
where
$\mathcal D:=\left\{\rho\in (0,1)^N 
:\sum_{i=1}^N\rho_i< 
1\right\}$.
We have
\begin{equation*}
\left(s'(\rho)\right)_\ell=\log\frac{\rho_\ell}{1-\sum_{j=1}^N\rho_j}
\quad \text{and }\quad
\left(s''(\rho)\right)_{k\ell}=\frac{\delta_{k\ell}}{\rho_\ell}+\frac{1}{1-\sum_{j=1}^N\rho_j}.
\end{equation*}
Then $s\in C^2(\mathcal D,[0,\infty))\cap
C^0(\overline{\mathcal D})$ and is
convex.
Moreover, $u:\mathbb R^N\to\mathcal D$ defined as
\begin{align*}
u_\ell(w)=\frac{e^{w_\ell}}{1+\sum_{i=1}^N e^{w_i}}\quad\mbox{for
}\ell=1,\ldots, N,
\end{align*}
{a choice first used in \cite{burger} to investigate the case $N=2$,}
is in $C^1(\mathbb R^N,\mathcal D)$, 
and is the inverse of $s'$.
Thus, the preamble of assumption~(H2) is satisfied.

In the numerical experiments below,
we use continuous space-time finite element discretization spaces.
On the space-time
cylinder  $Q_T=\Omega\times (0,T)$,
with $\Omega$ bounded interval ($d=1$) or Lipschitz polytope
($d>1$), 
we consider families of 
shape-regular
simplicial or Cartesian meshes $\{\mathcal T_h\}_{h>0}$.
The parameter
$h$ denotes the mesh granularity, namely
$\mathcal T_h=\{K_i, i=1,\ldots, N_h\}$, $h_K:=\text{diam}(K)$, and 
$ h:=\max_{K\in \mathcal T_h}h_{K}$.

As discretization spaces, we choose $\{\bm
V_h\}_{h>0}=\{\bm V_h^p,\ p\in\mathbb N\}_{h>0}$, with
\begin{align}
\label{eq.Vhp}
\bm V_h^p=\left\{v\in C^0(\overline Q_T)^N :\ v_{|_{K}}\in\mathcal P_p(K)^N\quad \forall K\in\mathcal T_h\right\},
\end{align}
where $\mathcal P_p(K)$ denotes the space of polynomial functions on
$K$ of degree at most $p$, if $K$ is a simplex, or of degree at most
$p$ in each variable, if $K$ is a cuboid. Therefore, the
approximability assumption (H4)
in the first part of section~\ref{sec:STGM} is satisfied.

Defining $B:\mathbb R^N\to \mathbb R^{N\times N}$ as
\begin{align*}
B(w)=A(u(w))u^{\prime}(w),
\end{align*}
the space-time Galerkin method~\eqref{eq.galerkin} can be rewritten
more explicitly in terms of the entropy variable unknown as follows:
\begin{multline}\label{MS.CG}
\text{Find}\ w_h^\epsilon\in \bm V_h^p \ \text{such that}\\
\epsilon(\phi,w_{h}^\epsilon)_{{\WeightedSpace}(Q_T)}+
\int_{\Omega}\phi(T) \cdot u(w_h^\epsilon(T))dx
    -\int_{\Omega}\phi(0) \cdot \rho_0 dx
    -\int_{Q_T}\partial_t\phi \cdot u(w_h^\epsilon) dxdt\\
+\sum_{i,j=1}^N\int_{Q_T}\nabla \phi_{i} \cdot
            B_{ij}(w_h^\epsilon)\nabla (w_h^\epsilon)_j dxdt
    =\int_{Q_T}\phi \cdot f(u(w_h^\epsilon)) dxdt\\
    \text{for all } \phi\in \bm V_h^p. 
\end{multline}
Throughout this section, we measure the absolute numerical error defined by
$\|\rho-u(w_h^\epsilon)\|_{L^2(Q_T)}$.

\begin{remark}[On the choice of $\varepsilon$]
{
    For the analysis the regularization term with the factor $\varepsilon>0$ is crucial, as it delivers essential bounds on the entropy variable $w$.
    In the numerical examples we will specify the choice of $\varepsilon$ for every example, and show on several occasions that it is possible to choose $\varepsilon=0$. 
    However, it is crucial to note that in these examples the entropy variable representing the exact solution has `nice' bounds.
  Conversely, when the solution approaches the singularities of the entropy, it is required to choose $\epsilon>0$ large enough for the solver to converge.
  In general, it is best to choose $\epsilon$ as small as possible, as we observe in the second example in \cref{sec:porousmedium}.
}
\end{remark}

\subsection{Heat equation}\label{sec:heat}
We apply our general approach to the linear heat equation:
\begin{equation*}
\begin{cases}
\partial_t \rho=\Delta\rho&\mbox{in }\Omega,\ t>0,\\
\partial_\nu \rho = 0&\mbox{on }\partial\Omega,\ t>0,\\
{\rho(0)=\rho_0}&{\mbox{in }\Omega.}
\end{cases}
\end{equation*}
This corresponds to problem~\eqref{eq.classical} with
$N=1$, $A\equiv1${, and $f\equiv0$}.
Furthermore,  $\mathcal D=(0,1)$ and the entropy density $s:\mathcal D\to [0,+\infty)$ is given by
\begin{equation*}
s(\rho)=\rho\log \rho+(1-\rho)\log(1-\rho)+\log(2),
\end{equation*}
and thus $s^{\prime}(\rho)=\log\frac{\rho}{1-\rho}$, and $s^{\prime\prime}(\rho)=\frac{1}{\rho(1-\rho)}$. 

For this choice of $A(\rho)$ and $f(\rho)$,
assumption (H1) is obviously satisfied, and
assumptions (H2a) and (H2b) are fulfilled with
$\gamma=4$ and $C_f=0$.

\begin{table}[ht]\centering
\caption{Numerical results for the heat equation.}
\label{tab:heateqresults}
\begin{tabular}{lll}
    \toprule  
    \multicolumn{3}{c}{$p=3$}\\
    $h$ & error & rate \\
    \midrule
    \begin{filecontents*}{heat2d.csv}
p,h,hnr,error,rate,steps
1,0.5,1,0.009778765622673427,0,6
1,0.25,2,0.01690442138165223,-0.7896763654955901,6
1,0.125,3,0.003916851066249056,2.1096344577049573,6
1,0.0625,4,0.0009258913037861054,2.0807795310360224,6
1,0.03125,5,0.00022779493859197662,2.0231071444088884,6
1,0.015625,6,5.671624503337852e-05,2.0059017664500027,6
2,0.5,1,0.009464016374323505,0,6
2,0.25,2,0.002438274704601351,1.9565918955715609,6
2,0.125,3,0.0002801517143913621,3.1215804479156133,6
2,0.0625,4,3.461704161901181e-05,3.016653979174455,6
2,0.03125,5,4.552277087495296e-06,2.926822157664495,6
2,0.015625,6,6.720218818896685e-07,2.760008259246122,6
3,0.5,1,0.0023012326630046727,0,6
3,0.25,2,0.0003077894753803451,2.902391046862833,6
3,0.125,3,2.411361700936186e-05,3.6740239261292276,6
3,0.0625,4,2.0995433419582784e-06,3.521700593550208,6
3,0.03125,5,1.3323390280741674e-07,3.9780424271873973,6
3,0.015625,6,8.355225352402789e-09,3.9951386878901594,6
4,0.5,1,0.0004800529366333892,0,6
4,0.25,2,2.9969991711882364e-05,4.001602914191829,6
4,0.125,3,2.3444086820072852e-06,3.6762226945993985,6
4,0.0625,4,5.878964670187316e-08,5.317518166937388,6
4,0.03125,5,1.8152389542034152e-09,5.017330728855576,6
4,0.015625,6,5.745840755247983e-11,4.981497654371164,6
\end{filecontents*}
\csvreader[head to column names,filter ifthen=\equal{\p}{3}]{heat2d.csv}{}
    {$2^{-\hnr}$ 
    & \num[round-precision=2,round-mode=figures, scientific-notation=true]{\error} 
    & \num[round-precision=2,round-mode=figures, scientific-notation=true]{\rate} 
    \tabularnewline}
    \\\addlinespace[-\normalbaselineskip]\bottomrule
\end{tabular}
\hskip10pt
\begin{tabular}{lll}
    \toprule  
    \multicolumn{3}{c}{$p=4$}\\
    $h$ & error & rate \\
    \midrule
    \csvreader[head to column names,filter ifthen=\equal{\p}{4}]{heat2d.csv}{}
    {$2^{-\hnr}$ 
    & \num[round-precision=2,round-mode=figures, scientific-notation=true]{\error} 
    & \num[round-precision=2,round-mode=figures, scientific-notation=true]{\rate} 
    \tabularnewline}
    \\\addlinespace[-\normalbaselineskip]\bottomrule
\end{tabular}
\end{table}

\begin{figure}[!ht]  
\centering
\resizebox{0.33\linewidth}{!}{
\begin{tikzpicture}
    \begin{semilogyaxis}
        [ xlabel={$p$},
        ylabel={Error},
        ymajorgrids=true,
        grid style=dashed,
        xtick={1,2,...,4},
        max space between ticks=20,
        cycle list name=paulcolors,
        legend style={at={(0.03,0.03)},anchor=south west}
        ]
        \foreach \hnr in {1,2,...,4}{
            \edef\temp{\noexpand\addlegendentry{$h$=$2^{-\hnr}$}};
            \addplot+[discard if not={hnr}{\hnr}] table [x=p, y=error, col sep=comma] {./heat2d.csv};
            \temp
        }
    \end{semilogyaxis}
\end{tikzpicture}
}
\resizebox{0.335\linewidth}{!}{
\begin{tikzpicture}
    \begin{loglogaxis}
        [ xlabel={$h$},
        ylabel={Error},
        ymajorgrids=true,
        grid style=dashed,
        x dir=reverse,
        xmin=0.01,xmax=1,
        max space between ticks=20,
        cycle list name=paulcolors,
        legend style={at={(0.03,0.03)},anchor=south west}
        ]
        \foreach \p in {1,2,...,4}{
            \edef\temp{\noexpand\addlegendentry{$p$=\p}};
            \addplot+[discard if not={p}{\p}] table [x=h, y=error, col sep=comma] {./heat2d.csv};
            \temp
        }
        \logLogSlopeTriangle{0.45}{0.11}{0.04}{5}{cyan}; 
        \logLogSlopeTriangle{0.45}{0.11}{0.04}{2}{teal}; 
    \end{loglogaxis}
\end{tikzpicture}
}
\resizebox{0.32\linewidth}{!}{
\begin{tikzpicture}
    \begin{semilogyaxis}
        [ xlabel={Time},
        ylabel={Entropy},
        ymajorgrids=true,
        cycle list name=paulcolors,
        grid style=dashed,
        ]
        \addplot table [mark=none,x=time, y=entropy, col sep=comma] {./heat2d_entropy.csv};
    \end{semilogyaxis}
\end{tikzpicture}
}

\caption{Convergence rates for the space-time Galerkin approximation towards the exact solution of the heat equation, in polynomial degree $p$ (left), and mesh size $h$ (middle). {On the right we are plotting the entropy on a logarithmic scale, showing exponential convergence.}}
\label{fig.heateqresults}
\end{figure}
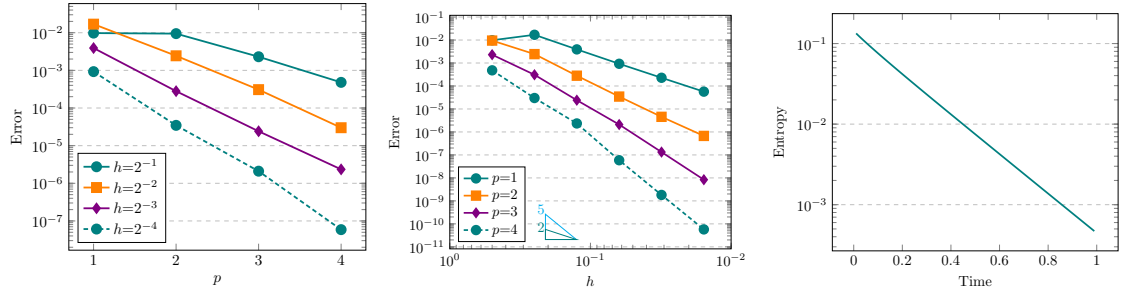

For the numerical tests, we take
$\Omega=(0,1)^2$ and
$\rho_0(\bm x)=0.5\cos(\pi x_1)\cos(\pi x_2)+0.5$, so that the problem has the
analytical solution given by
$$\rho(t,\bm x) = 0.5\exp(-8\pi^2t/\tau)\cos(2\pi x_1)\cos(2\pi x_2)+0.5,$$
where we use $\tau=7$ to rescale the time. 
The solution is shifted and scaled in order to avoid the singularities of $s'$ at 0 and 1.
Without this rescaling, the system matrix is highly ill-conditioned, which prohibits optimal convergence rates.
We solve \eqref{MS.CG}, setting $\epsilon=0$ and solving the
nonlinearity by Newton's method. We use
unstructured space-time simplicial meshes.
The Newton method converges in 6 steps, for all considered values of $h$ and $p$.
We measure the $L^2$ error on the whole space-time domain.
In \Cref{fig.heateqresults}, the convergence rates of the $h$- and the
$p$-version of the method are shown. 
We observe optimal rates, exponential in $p$ and of order $p+1$ in $h$.
In the case of $p=4$, we observe a preasymptotic region for very large mesh sizes; the exact rates are shown in \Cref{tab:heateqresults}.

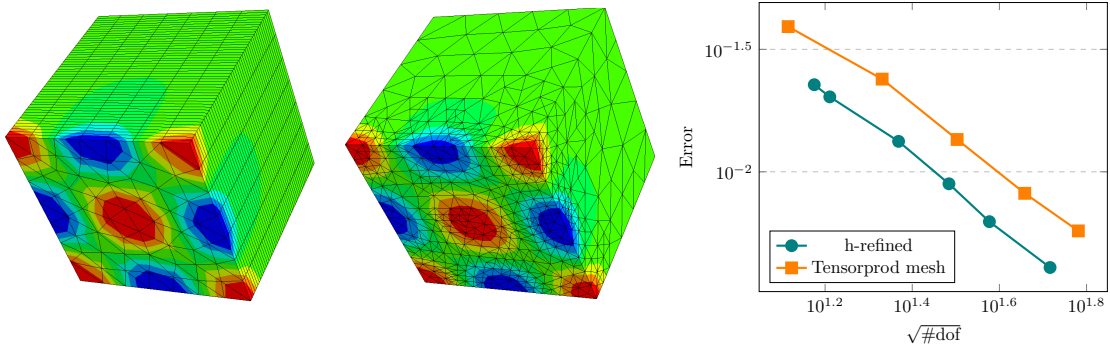
\begin{figure}[ht!]  
\centering
\includegraphics[trim=0 -50 50 100,clip,width=0.28\linewidth]{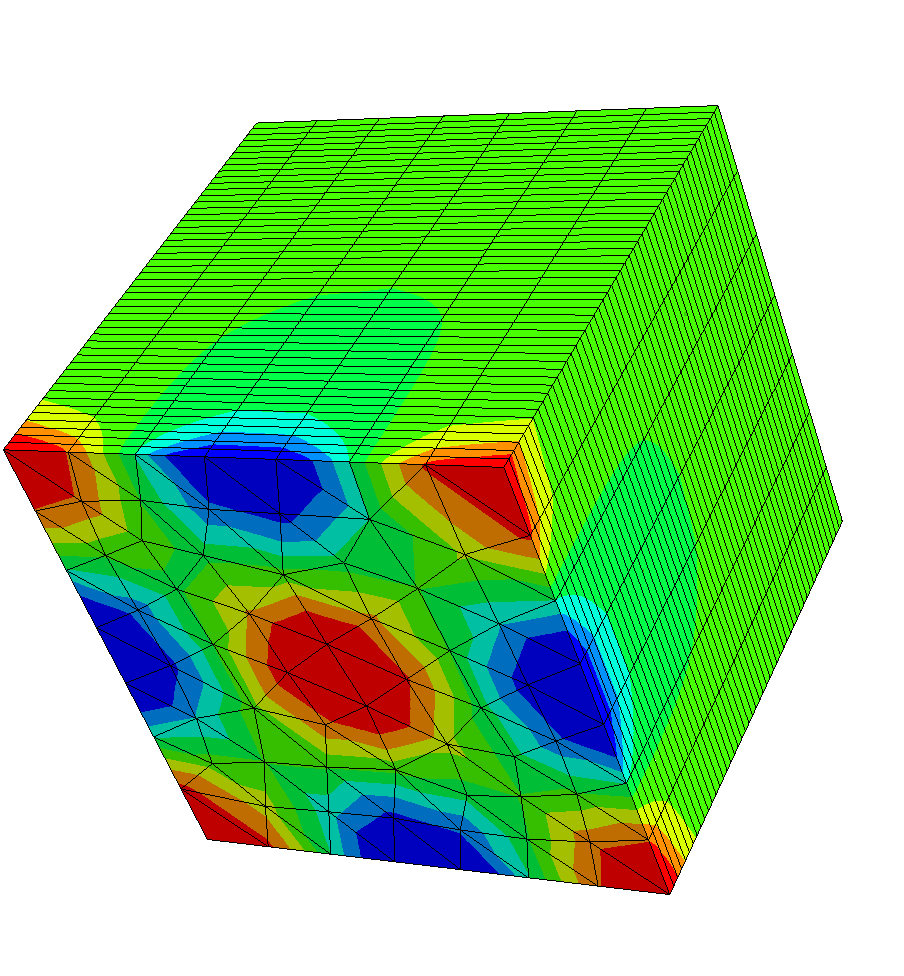}\quad
\centering
\includegraphics[trim=0 -50 50 100,clip,width=0.28\linewidth]{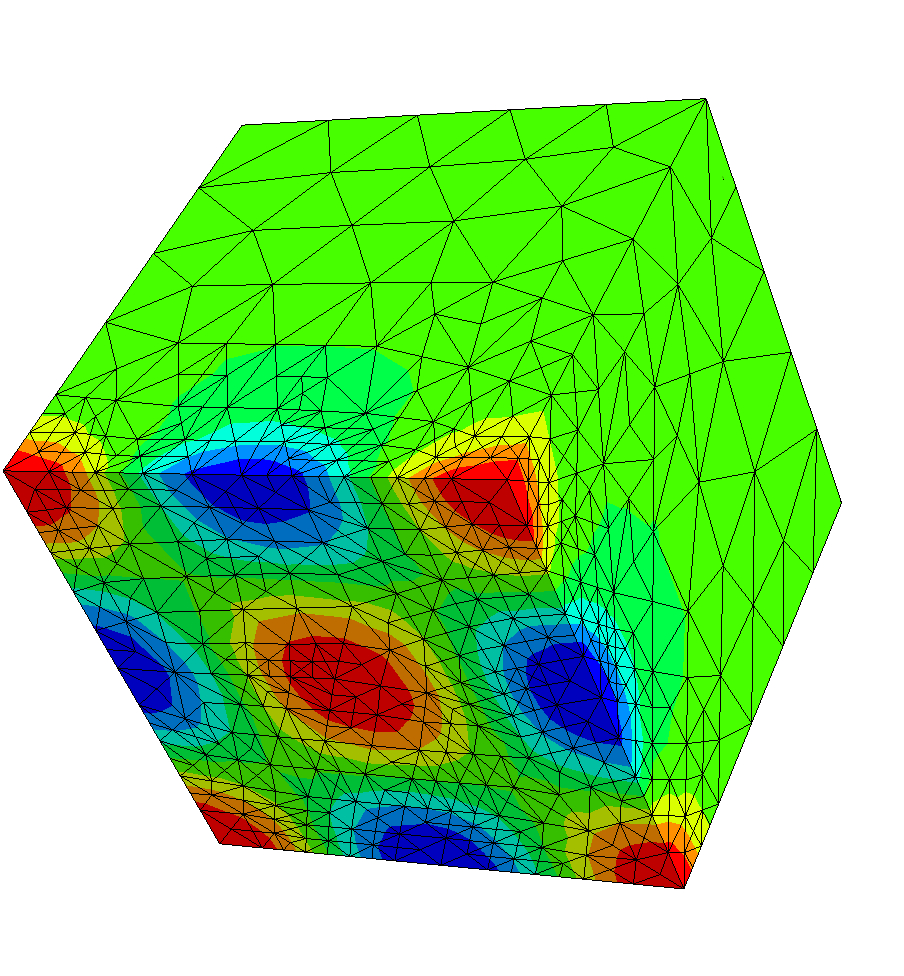}
\centering
\resizebox{0.4\linewidth}{!}{
\begin{tikzpicture}
    \begin{loglogaxis}
        [ xlabel={$\sqrt{\#\text{dof}}$},
        ylabel={Error},
        ymajorgrids=true,
        grid style=dashed,
        cycle list name=paulcolors,
        legend style={at={(0.03,0.03)},anchor=south west}
        ]
        \addplot+[] table [x=ndof, y=error, col sep=comma] {./aheat2d.csv};
        \addlegendentry{h-refined}
        \addplot+[] table [x=ndof, y=error, col sep=comma] {./aheat2dsq.csv};
        \addlegendentry{Tensorprod mesh}
    \end{loglogaxis}
\end{tikzpicture}
}
\caption{{Comparison of a mesh made from time slabs (left) and an
adapted space-time mesh (middle). The convergence of the two methods with respect to the number of degrees of freedom is shown on the right for $p=1$.}}
\label{fig.heateqadap}
\end{figure}

{In \Cref{fig.heateqadap} we highlight another feature of the
  space-time approach, namely the ability to use adapted mesh in space and time.
As a comparison, we use a mesh made of time slabs, a mesh structure similar to what would result from classic time-stepping methods. 
The timeslab height is given by $h_t\approx h_x^2$, with $h_x$ being the mesh size of the spatial mesh.
For the space-time adapted mesh, we start with a unstructured simplicial mesh of size $h=0.2$ and then apply adaptive refinement. 
For this simple example, we use a flux based error estimator and D\"orfler marking with $\theta=0.5$. 
We observe the same rate of convergence on both meshes, however using
the space-time adapted mesh allows us to obtain a given accuracy with fewer degrees of freedom.}

\subsection{The porous medium equation}\label{sec:porousmedium}
Let $m>1$.
The porous medium equation is given by
\begin{equation*}
\begin{cases}
\partial_t \rho=\Delta\rho^m&\mbox{in }\Omega\ t>0,\\
\partial_\nu (\rho^m) = 0&\mbox{on }\partial\Omega,\ t>0,\\
    \rho(0)=\rho_0&\mbox{in }\Omega.
\end{cases}
\end{equation*}
We can write {it} in the form of
\eqref{eq.classical} for {$N=1$}, $A(\rho)=m\rho^{m-1}$, and $f\equiv0$.
The entropy density is the same as for the heat equation. 
\begin{proposition}
Assumptions (H1) and (H2) are satisfied for $m\in(1,2]$.
\end{proposition}
\begin{proof}
For $\mathcal D=(0,1)$ and $m>1$, $A(\rho)=m\rho^{m-1}$ is in $C^0(\overline{\mathcal D})$, thus (H1) is stisfied. As~(H2b) is obvious,  we only
{need} to prove that~(H2a) is satisfied, namely
that $s''(\rho)A(\rho)\geq \gamma$ for some $\gamma>0$ and all
$\rho\in\mathcal D$. Thus let $\rho\in(0,1)=\mathcal D$. Then,
whenever $m\in (1,2]$,
\begin{align*}
s''(\rho)A(\rho)&=\frac{m \rho^{m-1}}{\rho(1-\rho)}=\frac{m}{\rho^{2-m}(1-\rho)}\geq m=:\gamma. 
\end{align*}
\end{proof}

We test the space-time Galerkin method for this problem
with initial conditions and Neumann boundary conditions chosen such that
\begin{equation*}
\rho(x,t)=\left[\frac{(m-1)(x-\alpha)^2}{2m(m+1)(\beta-t)}\right]^{\frac{1}{m-1}}
\end{equation*}
is the exact solution, with $\alpha$ and $\beta$ real parameters, on $\Omega=(0,1)$.
We consider the case $m=2,\ \alpha=2,\ \beta=5$
on unstructured simplicial space-time meshes.

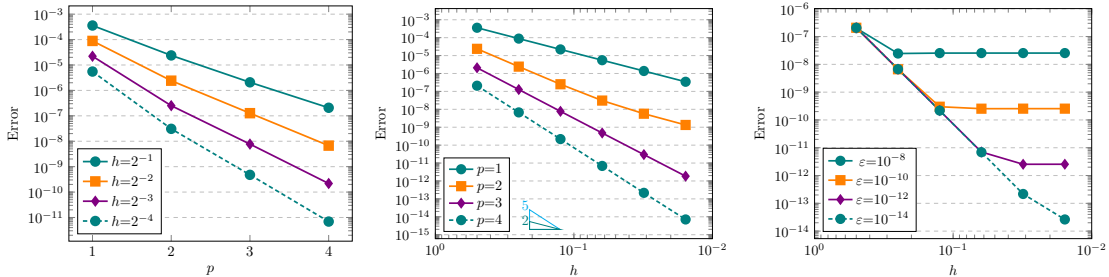
\begin{figure}[!ht]  
\centering
\resizebox{0.32\linewidth}{!}{
\begin{tikzpicture}
    \begin{semilogyaxis}
    [ xlabel={$p$}, ylabel={Error}, ymajorgrids=true, grid style=dashed, xtick={1,2,...,4}, max space between ticks=20,
        cycle list name=paulcolors,
        legend style={at={(0.03,0.03)},anchor=south west}
    ]
        \foreach \hnr in {1,2,...,4}{
            \edef\temp{\noexpand\addlegendentry{$h$=$2^{-\hnr}$}};
            \addplot+[discard if not={hnr}{\hnr}] table [x=p, y=error0, col sep=comma] {./porousm.csv};
            \temp
        }
    \end{semilogyaxis}
\end{tikzpicture}
}
\resizebox{0.33\linewidth}{!}{
\begin{tikzpicture}
    \begin{loglogaxis}
        [ 
        xlabel={$h$}, ylabel={Error}, ymajorgrids=true, grid style=dashed, x dir=reverse ,
        xmin=0.01,xmax=1,
        max space between ticks=20,
        cycle list name=paulcolors,
        legend style={at={(0.03,0.03)},anchor=south west}
        ]
        \foreach \p in {1,2,...,4}{
            \edef\temp{\noexpand\addlegendentry{$p$={\p}}};
            \addplot+[discard if not={p}{\p}] table [x=h, y=error0, col sep=comma] {./porousm.csv};
            \temp
        }
        \logLogSlopeTriangle{0.45}{0.11}{0.04}{5}{cyan}; 
        \logLogSlopeTriangle{0.45}{0.11}{0.04}{2}{teal}; 
    \end{loglogaxis}
\end{tikzpicture}
}
\resizebox{0.33\linewidth}{!}{
\begin{tikzpicture}
    \begin{loglogaxis}
        [ 
        xlabel={$h$}, ylabel={Error}, ymajorgrids=true, grid style=dashed, x dir=reverse ,
        xmin=0.01,xmax=1,
        max space between ticks=20,
        cycle list name=paulcolors,
        legend style={at={(0.03,0.03)},anchor=south west}
        ]
        \addlegendentry{$\epsilon$=$10^{-8}$};
        \addplot+[discard if not={p}{4}] table [x=h, y=error8, col sep=comma] {./porousm.csv};
        \addlegendentry{$\epsilon$=$10^{-10}$};
        \addplot+[discard if not={p}{4}] table [x=h, y=error10, col sep=comma] {./porousm.csv};
        \addlegendentry{$\epsilon$=$10^{-12}$};
        \addplot+[discard if not={p}{4}] table [x=h, y=error12, col sep=comma] {./porousm.csv};
        \addlegendentry{$\epsilon$=$10^{-14}$};
        \addplot+[discard if not={p}{4}] table [x=h, y=error14, col sep=comma] {./porousm.csv};
    \end{loglogaxis}
\end{tikzpicture}
}
\caption{{Convergence rates towards the exact solution of the
    porous medium equation. Convergence in terms of the polynomial
    degree $p$ for different mesh sizes is shown on the left. We
    consider convergence in mesh size $h$ for different values of the polynomial order $p$ and fixed $\epsilon=0$ in the middle, and for fixed $p=4$ and different values of regularization parameter $\epsilon$ on the right.}}
\label{fig:porousmexactsol}
\end{figure}

In \Cref{fig:porousmexactsol}, we show the convergence rates of the scheme.
Regardless of the nonlinearity, we match the convergence rates of the heat equation, i.e. exponential in $p$ and of order $p+1$ in $h$.
{The convergence rates in terms of $h$ are also considered for different values of $\epsilon$. 
    We observe that $\epsilon$ introduces a lower bound on the error. 
    Therefore, choosing it as small as possible, such that the solver
    still converges, gives the best results.
On the other hand, in the next example, we can see that for certain solutions, that produce a very ill-conditioned system, we must choose $\epsilon$ fairly large.
}
    
\begin{figure}[!ht]  
\centering
\resizebox{0.4\linewidth}{!}{
    \begin{tikzpicture}
    \begin{axis}
    [ xlabel={$x$},
    ylabel={$\rho$},
    ymajorgrids=true,
    grid style=dashed,
    yticklabel style={ /pgf/number format/precision=3 },
    cycle list name=paulcolors
    ]
        \foreach \t in {0,0.042,0.083,0.167}{
            \addplot+ [discard if not={t}{\t}] table [mark=none, x=x, y=val, col sep=comma] {./porousmwaitingtime.csv};
        }
        \addplot [mark=none] coordinates { (0.785398,0) (0.785398,1)};
        \addplot [mark=none] coordinates { (3.9269906,0) (3.9269906,1)};
    \end{axis}
    \end{tikzpicture}
}
\resizebox{0.42\linewidth}{!}{
    \begin{tikzpicture}
    \begin{axis}
    [ xlabel={Time},
    ylabel={Value at interface},
    ymajorgrids=true,
    grid style=dashed,
    yticklabels={-0.04,-0.02,0,0.02,0.04,0.06,0.08,0.1,0.12},
    xticklabel style={
        /pgf/number format/fixed,
        /pgf/number format/precision=2,
    },
    ]
    \addplot [line width=1pt,color=teal,discard if not={x}{0.785}] table [ mark=none, x=t, y=val, col sep=comma] {./porousmwaitingtime.csv};
    \addplot [color=red,mark=none] coordinates { (0.083333,-0.01) (0.083333,0.01)};
    \end{axis}
    \end{tikzpicture}
}
\caption{Snapshots of the solution of the porous medium equation emitting a waiting time, at different times (left) and the value at the left interface (right).}
\label{fig.waitingtime}
\end{figure}
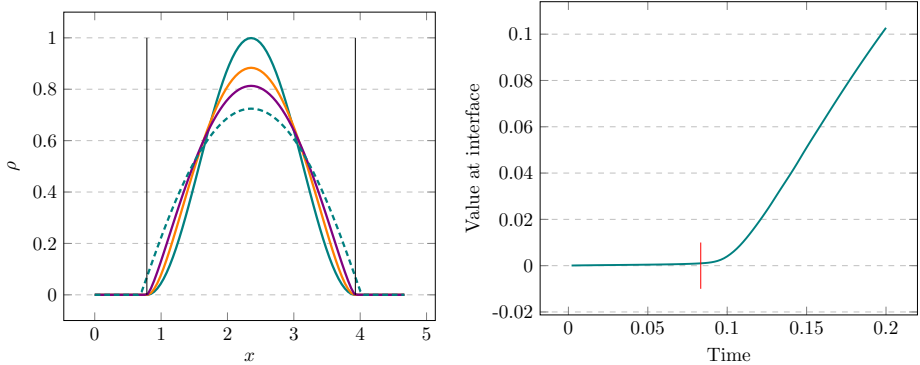

In contrast to the heat equation, the power law in the porous medium equation introduces a finite propagation speed of the solution.
This is best observed by the interesting behavior of certain initial
conditions that induce a waiting time.
That is, the solution keeps a fixed support until the waiting time is reached.
On $\Omega=(0,\pi)$, the initial condition given by 
\begin{align*}
\rho_0(x)=
\begin{cases}
    \sin^{2/(m-1)}(x) & \text{ if } 0\leq x\leq\pi,\\
    0 & \text{ otherwise, }
\end{cases}
\end{align*}
produces this behavior.
It is shown in \cite{MR2019614} that the corresponding solution has a waiting time of $t^*=\frac{m-1}{2m(m+1)}$.
As we choose $m=2$, here $t^*=0.08\dot{3}$.
We {modify the initial condition to $\rho_0(x)=10^{-16}$ for $x\notin [0,\pi]$ to avoid ill-conditioning.}
Furthermore, to ensure convergence of the Newton method used as a
nonlinear solver, we had to choose $\epsilon = {10^{-6}}$, making use of the regularization term.
We solve on a Cartesian space-time mesh until final time $T=0.2$, with spatial mesh size $h_s=0.05$, and temporal mesh size $h_t=h_s/2$, and fix $p=5$.
The results are shown in \Cref{fig.waitingtime}.
Looking at snapshots of the numerical solution we can observe that it keeps a compact support set. 
In \Cref{fig.waitingtime}, on the right, we plot the value of the
solution on the left interface against time, marking the expected
waiting time $t^*$ with the vertical line.

\subsection{The Fisher-KPP equation}\label{sec:FisherKPP}
We consider the Fisher-KPP equation
\begin{equation*}
\begin{cases}
\partial_t \rho=A\Delta\rho+\rho(1-\rho)&\mbox{in }\Omega,\ t>0,\\
A\partial_\nu \rho = 0&\mbox{on }\partial\Omega,\ t>0,\\
\rho(0)=\rho_0&\mbox{in }\Omega,
\end{cases}
\end{equation*}
with $A>0$ now constant. 
This agrees with formulation~\eqref{eq.classical}, with
$N=1$, $A(\rho)=A$, and $f(\rho)=\rho(1-\rho)$. We set again
$\mathcal D:=(0,1)$. 
Assumptions~(H1) and (H2a) are clearly satisfied.
Choosing an entropy density such that assumption~(H2b) is
satisfied with $C_f=0$ allows for the right-hand side of the entropy estimate \eqref{w.entropy} to be independent of time. 
Motivated by this, we now investigate the rescaled entropy density
$s:\mathcal D\to (0,+\infty)$ 
given by
\begin{equation}\label{eq.fischerentropy}
s(\rho)={\rho}\log {\rho}+(n-{\rho})\log(n-{\rho}),
\end{equation}
with $n$ to be chosen.
Note that 
$f(\rho)> 0 $ for $\rho\in (0,1)$,
and $n/\rho-1>1$ if and only if $\rho<n/2$.
Thus, 
\begin{align*}
f(\rho) s'(\rho)&=\rho(1-\rho)\log\frac{\rho}{n-\rho}
              =-\rho (1-\rho) 
              \log{(\frac n \rho-1)}\leq 0
\end{align*}
for all for $\rho\in(0,1)$ if and only if $n\geq 2$. 
We choose $n=2$ so that the hypothesis (H2b) is fulfilled with $C_f=0$.

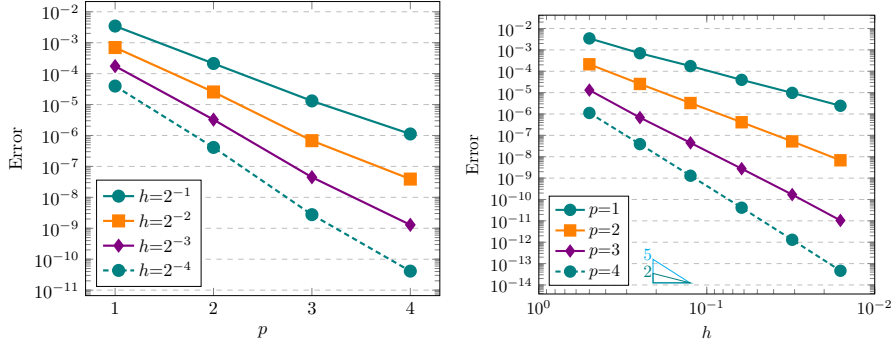
\begin{figure}[!ht]  
\centering
\resizebox{0.4\linewidth}{!}{
\begin{tikzpicture}
    \begin{semilogyaxis}
        [ xlabel={$p$}, ylabel={Error}, ymajorgrids=true, grid style=dashed, xtick={1,2,...,4}, 
        max space between ticks=20,
        cycle list name=paulcolors,
        legend style={at={(0.03,0.03)},anchor=south west}
        ]
        \foreach \hnr in {1,2,...,4}{
            \edef\temp{\noexpand\addlegendentry{$h$=$2^{-\hnr}$}};
            \addplot+[discard if not={hnr}{\hnr}] table [x=p, y=error, col sep=comma] {./fisherkpp.csv};
            \temp
        }
    \end{semilogyaxis}
\end{tikzpicture}
}
\resizebox{0.4\linewidth}{!}{
\begin{tikzpicture}
    \begin{loglogaxis}
        [ xlabel={$h$}, ylabel={Error}, ymajorgrids=true, grid style=dashed, x dir=reverse, 
        xmode=log,
        ymode=log,
        xmin=0.01,xmax=1,
        max space between ticks=20,
        cycle list name=paulcolors,
        legend style={at={(0.03,0.03)},anchor=south west}
        ]
        \foreach \p in {1,2,...,4}{
            \edef\temp{\noexpand\addlegendentry{$p$=\p}};
            \addplot+[discard if not={p}{\p}] table [x=h, y=error, col sep=comma] {./fisherkpp.csv};
            \temp
        }
        \logLogSlopeTriangle{0.45}{0.11}{0.04}{5}{cyan}; 
        \logLogSlopeTriangle{0.45}{0.11}{0.04}{2}{teal}; 
    \end{loglogaxis}
\end{tikzpicture}
}
\caption{Convergence rates in polynomial degree $p$ (left) and mesh size $h$ for the exact solution of the Fisher-KPP equation.}
\label{fig:fisherkppexactsol}
\end{figure}

We start again by investigate convergence towards a smooth solution. 
We choose $\Omega=(0,1)$, and initial conditions and
Neumann boundary conditions such that 
\begin{equation*}
\rho(x,t)= \frac{1}{\left[1+\exp(-\frac56 t+\frac{1}{\sqrt{6}}x)\right]^2}
\end{equation*}
is the exact solution for $A=1$.
We set {$\epsilon=0$} and solve on unstructured simplicial space-time meshes.
The results are presented in \Cref{fig:fisherkppexactsol}.
We observe again optimal convergence rates in both $p$ and $h$, namely
exponential in $p$ and of order $p+1$ in $h$.

Next, we aim to reproduce the experiments presented in \cite{bonizzoni2019structurepreserving},
considering an initial condition with a jump, given by $\rho_0(x)= 1 \text{ if } 0<x<1/2$ and 0 elsewhere, with diffusion coefficient $A=10^{-4}$.
We solve using $p=3$ on a Cartesian mesh with $h_s=0.025$, $h_t=0.4$ up to $T=8$.
We choose $\epsilon = 10^{-8}$ to avoid ill-conditioning in the solver.

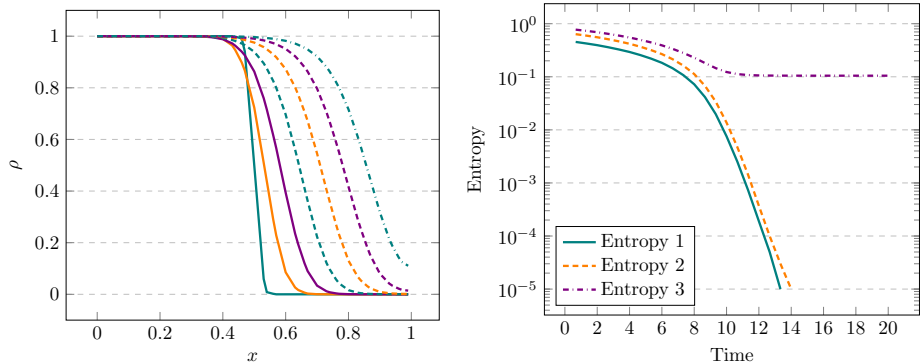
\begin{figure}[!ht]  
\centering
\resizebox{0.4\linewidth}{!}{
    \begin{tikzpicture}
    \begin{axis}
    [ xlabel={$x$},
    ylabel={$\rho$},
    ymajorgrids=true,
    grid style=dashed,
    yticklabel style={
        /pgf/number format/precision=3,
    },
    cycle list name=paulcolors
    ]
        \foreach \t in {0.0,1.333,2.667,4.0,5.333,6.667,8.0}{
            \addplot+ [discard if not={t}{\t}] table [mark=none, x=x, y=val, col sep=comma] {./fisherkppspecies.csv};
        }
    \end{axis}
    \end{tikzpicture}
}
\resizebox{0.42\linewidth}{!}{
    \begin{tikzpicture}
    \begin{semilogyaxis}
    [ xlabel={Time},
    ylabel={Entropy},
    ymajorgrids=true,
    grid style=dashed,
    max space between ticks=20,
    legend pos=south west,
    cycle list name=paulcolorsnm,
    ]
    \foreach \n in {0,1,2}{
        \addplot+ [discard if not={n}{\n}] table [mark=none, x=t, y=entropy, col sep=comma] {./fisherkppentropy.csv};
    }
    \addlegendentry{Entropy 1}
    \addlegendentry{Entropy 2}
    \addlegendentry{Entropy 3}
    \end{semilogyaxis}
    \end{tikzpicture}
}
\caption{Snapshots of the numerical solution for the Fisher-KPP (left) and different choices of the entropy (right). The choices are as follows: Entropy 1 is the one used in \cite{bonizzoni2019structurepreserving}, Entropy 2 is given by \eqref{eq.fischerentropy} with $n=2$, and Entropy 3 is \eqref{eq.fischerentropy} with $n=2.1$.}
\label{fig:fisherkppspecies}
\end{figure}

Snapshots of the numerical solution are taken every $1.\dot{3}$ seconds, the results are shown in \Cref{fig:fisherkppspecies} on the left.
In \Cref{fig:fisherkppspecies} on the right, we consider different choices for the entropy up to $T=15$.
Note that at the point in time the solution has already converged to $\rho\equiv 1$.
The choice for the entropy density in~\cite{bonizzoni2019structurepreserving} was $\rho\log(\rho)-\rho+1$.
{This choice is not covered by our assumptions, however, using it produces the correct results, as conjectured in \Cref{rem:as}.}
We compare it to the entropy in \eqref{eq.fischerentropy} for different values of $n$ in \Cref{fig:fisherkppspecies}.
For the choice of $n=2$, we recover a similar behavior of the
entropy, namely, a region with slow decay followed by an exponential
decay.
As the solution converges to $1$ it can easily be seen that for $n>2$ the entropy does not {converge} to zero exponentially, as exemplified by the choice of $n=2.1$ in the figure, {since it is not the correct relative entropy with respect to the equilibrium.}

\subsection{{The three-component Maxwell-Stefan system}}\label{sec.ThreeGas}
{The Maxwell-Stefan system for $N=2$} can be written as
\begin{equation*}
\begin{cases}
\partial_t\rho_i=\nabla\cdot \left(\sum_{j=1}^{2}A_{ij}(\rho_1,\rho_2)\nabla\rho_{j}\right)&\mbox{in }\Omega,\ t>0,\\
\sum_{j=1}^2 A_{ij}(\rho_1,\rho_2)\partial_\nu \rho_j = 0&\mbox{on }\partial\Omega,\ t>0,\\
\rho_i(0)=(\rho_0)_i&\mbox{in }\Omega
\end{cases}
\end{equation*}
for $i=1,2$, with
\begin{equation}
A(\rho_1,\rho_2)=\frac{1}{\delta(\rho_1,\rho_2)}\begin{pmatrix}
d_1+(d_3-d_1)\rho_1 &(d_3-d_2)\rho_1\\
(d_3-d_1)\rho_2 & d_2+(d_3-d_2)\rho_2
\end{pmatrix}
\end{equation}
and
\[\delta(\rho_1,\rho_2)=
d_1d_2(1-\rho_1-\rho_2)+d_2d_3\rho_1+d_3d_1\rho_2. \]
The unknowns $\rho_1$ and $\rho_2$ represent the concentrations of
the first two gases ($\rho_3=1-(\rho_1+\rho_2)$); the parameters
$d_1$, $d_2$, and $d_3$ are the related to the binary diffusion coefficients of the three gases.
In section~\ref{sec.MaxwellStefan} below, we derive this form of the Maxwell-Stefan system, prove that it fits our framework, and discuss the case $N>2$.

\begin{figure}[!ht]  
\centering
\includegraphics[trim=0 600 0 50,clip,width=0.8\textwidth]{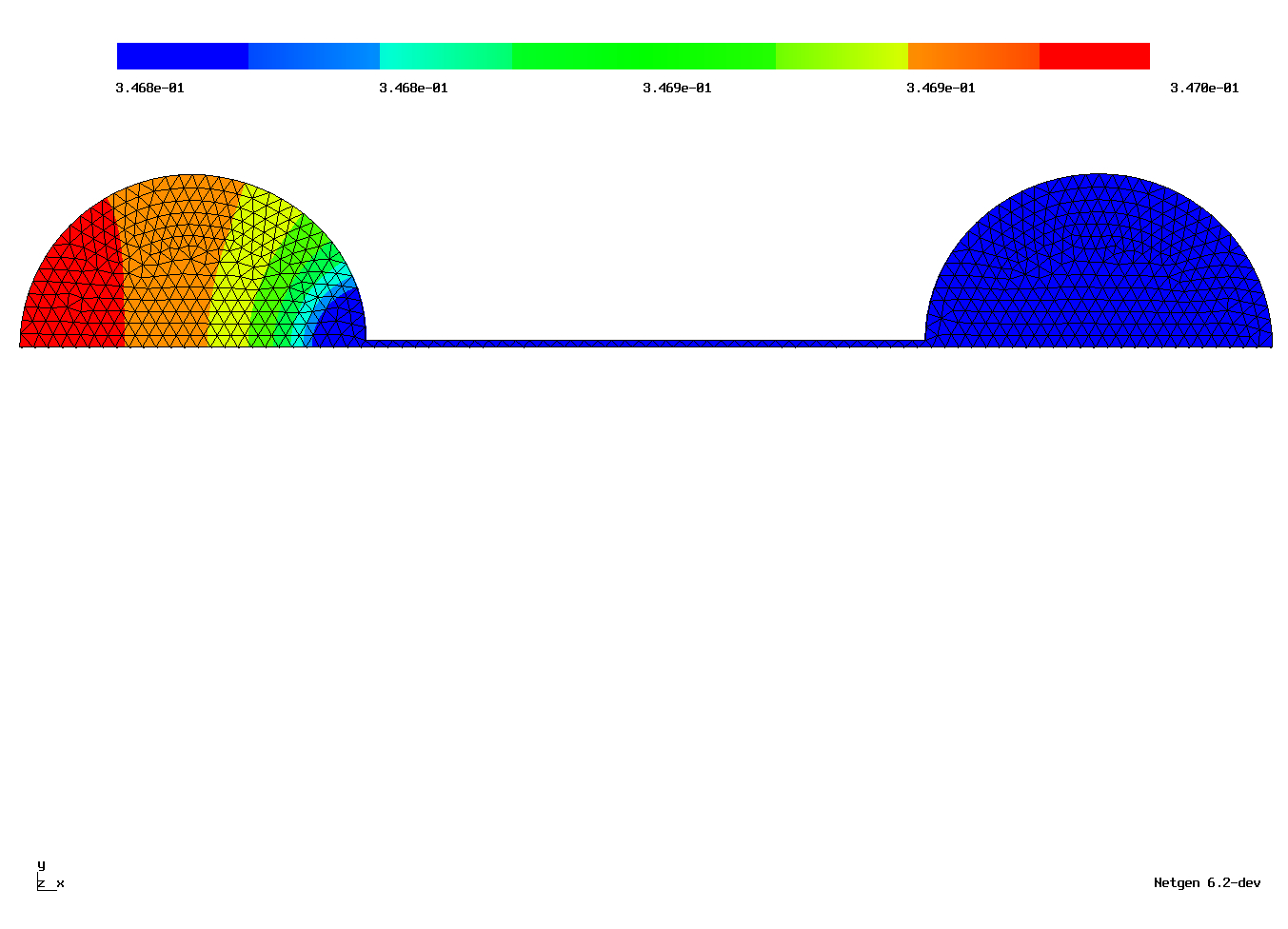}
\caption{The mesh used for the  Duncan-Toor example, depicting the carbon dioxide content after about ten hours.}
\label{fig.dunctoormesh}
\end{figure}

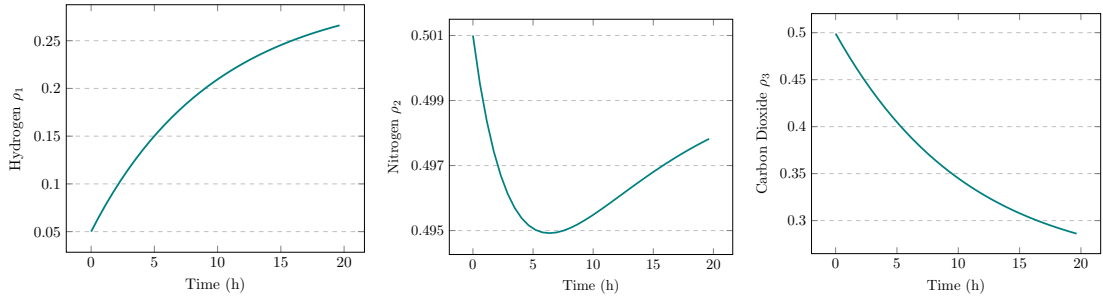
\begin{figure}[!ht]  
\centering
\resizebox{0.33\linewidth}{!}{
    \begin{tikzpicture}
    \begin{axis}
    [ xlabel={Time (h)},
    ylabel={Hydrogen $\rho_1$},
    ymajorgrids=true,
    cycle list name=paulcolors,
    yticklabels={0,0.05,0.1,0.15,0.2,0.25},
    grid style=dashed]
    \addplot table [mark=none,x=time, y=hydro, col sep=comma] {./duncantoorbig.csv};
    \end{axis}
    \end{tikzpicture}
}
\resizebox{0.32\linewidth}{!}{
    \begin{tikzpicture}
    \begin{axis}
    [ xlabel={Time (h)},
    ylabel={Nitrogen $\rho_2$},
    ymajorgrids=true,
    grid style=dashed,
    cycle list name=paulcolors,
    yticklabel style={
        /pgf/number format/precision=3,
        /pgfplots/ytick={0.495,0.497,...,0.501},
    },
    ymax=0.5015
    ]
    \addplot table [mark=none, x=time, y=nitro, col sep=comma] {./duncantoorbig.csv};
    \end{axis}
    \end{tikzpicture}
}
\resizebox{0.32\linewidth}{!}{
    \begin{tikzpicture}
    \begin{axis}
    [ xlabel={Time (h)},
    ylabel={Carbon Dioxide $\rho_3$},
    ymajorgrids=true,
    cycle list name=paulcolors,
    grid style=dashed]
    \addplot table [mark=none, x=time, y=carbo, col sep=comma] {./duncantoorbig.csv};
    \end{axis}
    \end{tikzpicture}
}
\caption{Comparison of the mole fractions in the left side of the device.}
\label{fig.dunctoor}
\end{figure}

In \cite[Sec. 2]{boudin} numerical results were presented for the three component gas diffusion experiment originally performed by Duncan and Toor in \cite{DuncanToor}.
The setting is the following. 
Consider two spherical bulbs of volume 77.99 cm$^3$ (radius 26.49
  mm) and 78.63 cm$^3$ (radius 26.58 mm), respectively, 
which are connected by a capillary tube of length 85.9 mm and diameter 2.08 mm, with a valve in the middle.
We consider the Maxwell-Stefan equations with $N=2$, corresponding to the gas mixture composed of hydrogen ($\rho_1$), nitrogen ($\rho_2$), and carbon dioxide ($\rho_3$). 
We consider the following initial gas mixture in the left- and right-hand side of the device:
\begin{align*}
\text{Left:}\quad (\rho_0)_1=0.000,\quad (\rho_0)_2=0.501,\quad (\rho_0)_3=0.499,\\
\text{Right:}\quad (\rho_0)_1=0.501,\quad (\rho_0)_2=0.499,\quad (\rho_0)_3=0.000.
\end{align*}
For these gases, the coefficients $d_1$, $d_2$, and $d_3$ are
  given in terms of the binary diffusion coefficients (see
  section~\ref{subsec.maxstef} below) as follows:
%
%
$$
{
d_{1}^{-1}=D_{13}={68.0} \text{ mm}^2\text{s}^{-1}, \quad
d_{2}^{-1}=D_{23}={16.8} \text{ mm}^2\text{s}^{-1}, \quad
d_{3}^{-1}=D_{12}={83.3} \text{ mm}^2\text{s}^{-1}.}$$

{As in \cite{boudin}, we can reduce the domain to two dimensions, as the device and initial conditions are axially symmetric and the flux vector has no angular component.} In \Cref{fig.dunctoormesh}, the computational domain is shown. 
We choose the spatial mesh size $h_s=2.08$ mm, equal to the diameter of the tube. 
The size of the Cartesian product mesh in time is chosen as 20.8 s.
We solve iteratively on these slabs, restarting the computations with the previous solution as initial condition.
{We fix $p=2$ and $\epsilon=10^{-10}$}.

The results are shown  in \Cref{fig.dunctoor}.
We recover the same behavior shown in \cite{boudin}.
Both hydrogen and carbon dioxide converge monotonically to the expected equilibrium. 
Nitrogen shows the peculiar behavior known from the experiment. 

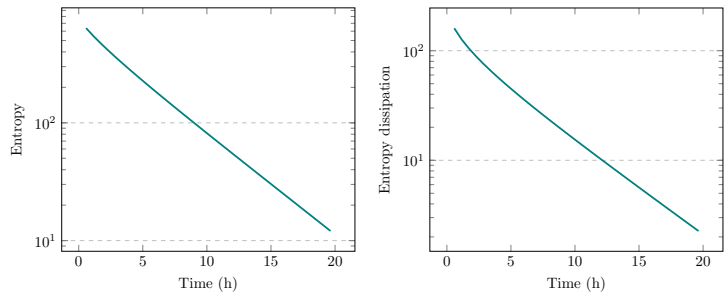
\begin{figure}[!ht]  
\centering
\resizebox{0.32\linewidth}{!}{
    \begin{tikzpicture}
    \begin{semilogyaxis}
    [ xlabel={Time (h)},
    ylabel={Entropy},
    ymajorgrids=true,
    cycle list name=paulcolors,
    grid style=dashed]
    \addplot table [mark=none,x=time, y=entropy, col sep=comma] {./duncantoorbig.csv};
    \end{semilogyaxis}
    \end{tikzpicture}
}
\resizebox{0.32\linewidth}{!}{
    \begin{tikzpicture}
    \begin{semilogyaxis}
    [ xlabel={Time (h)},
    ylabel={Entropy dissipation},
    ymajorgrids=true,
    grid style=dashed,
    cycle list name=paulcolors,
    yticklabel style={
        /pgf/number format/precision=3,
    },
    ]
    \addplot table [mark=none, x=time, y=dtentropy, col sep=comma] {./duncantoorbig.csv};
    \end{semilogyaxis}
    \end{tikzpicture}
}
\caption{{Entropy and entropy dissipation (in absolutes) for the Duncan-Toor example, both showing exponential convergence to the equilibrium.}}
\label{fig.dunctoor_entropy}
\end{figure}

{In \Cref{fig.dunctoor_entropy}, we show the relative entropy and its dissipation, i.e. the time derivative of the entropy, both converge exponentially. 
This can be expected, as this behavior has been proven for similar types of equations, see for example \cite[eq. (1.16)]{MR1842428} where it is formally stated that the entropy dissipation converges exponentially for the Fokker-Planck equation by the use of the Bakry-Emery method.
}

\section{The Maxwell-Stefan system revisited}\label{sec.MaxwellStefan}
In this section, we derive the formulation of the Maxwell-Stefan
system as that used in section~\ref{sec.ThreeGas}, and show that it
fits into the general framework of
section~\ref{sec:generalsetting} (section~\ref{subsec.maxstef}). For the case $N>2$,
in which an explicit representation of the currents may not be easily
derived, we introduce and analyze an alternative space-time Galerkin
method, which is based on a formulation that is implicit for the
currents (section~\ref{sec:implicit}).

Let $\rho_0\in L^\infty(\Omega)^{N+1}$ such that $ \rho_0\geq0$ and
$\sum_{i=1}^{N+1}  
(\rho_0)_i=1$. The Maxwell-Stefan equations are given by the continuity equations
\begin{equation}\label{continuity.eq}
\left\{
\begin{aligned}
\partial_t\rho_i +\nabla\cdot J_i = 0&\qquad \mbox{in } (0,T)\times\Omega,\\
\nu \cdot J_i =0 &\qquad\mbox{on } (0,T)\times\partial\Omega,\\
\rho_i(0)=(\rho_0)_i &\qquad\mbox{in }\Omega
\end{aligned}
\right.
\end{equation}
for $i=1,\ldots,N+1$,
where the currents $J_i$ are implicitly given by
\begin{align}\label{flux}
\nabla \rho_i = \sum_{j=1}^{N+1}
\frac{\rho_iJ_j-\rho_jJ_i}{D_{ij}}
\end{align}
for some $D_{ij}=D_{ji}>0$. 

\subsection{Explicit formula for the currents}\label{subsec.maxstef}
In this section, we establish an explicit representation of the
currents, which allows us to derive the formulation of the
Maxwell-Stefan system in the concentration variable unknowns.
We follow~\cite{Bothe2011} (see
also~\cite{JS13}).

Let $M _{ij}(\rho):=D_{ij}^{-1}\rho_i-\delta_{ij}
\sum_{k=1}^{N+1} 
D_{ik}^{-1}\rho_k$,  $i,j=1,\ldots,N+1$. Thus,
\[ \nabla \rho_i = \sum_{j=1}^{N+1}
M_{ij}(\rho)J_j.\]
Using $\rho_i\geq0$
and $D_{ij}=D_{ji}>0$, it is easy to see that $M(\rho)$ is quasi-positive
($M _{ij}(\rho)\ge 0$ for $i\ne j$).
Moreover, provided that $\rho_i>0$ for all $1\le i\le N+1$, $M(\rho)$
is irreducible.
Direct calculations show that 
\begin{equation*}
\mathrm{Ker}\,M(\rho)\supseteq\mathrm{span}\{\rho\} \qquad\mbox{and}\qquad
\mathrm{Im}\,M(\rho)\subseteq\left\{v: \sum_{i=1}^{N+1}
v_i=0\right\}.
\end{equation*}
Moreover, $R^{-1}M(\rho)R$, with $R=\text{diag}(\rho_1^{1/2},\ldots,\rho_{N+1}^{1/2})$, is symmetric, thus all the eigenvalues 
of $M(\rho)$ are real.
By the Perron-Frobenius theory for quasi-positive, irreducible
matrices, one deduces that the eigenvalue zero has multiplicity one
(we refer to~\cite{Bothe2011} or~\cite{JS13} for details). We deduce
\begin{equation}\label{ker.im.M}
\mathrm{Ker}\,M(\rho)=\mathrm{span}\{\rho\} \qquad\mbox{and}\qquad
\mathrm{Im}\,M(\rho)=\left\{v: \sum_{i=1}^{N+1}
v_i=0\right\}.
\end{equation}

As $M(\rho)$ is not invertible, we have to restrict ourselves to a subspace of all possible currents $J$
in order to obtain an explicit formula for $J$. For this, we make the
assumption that the total current
\[J_{\mathrm{tot}}:=\sum_{i=1}^{N+1}
J_i\]
vanishes. Then by summing in \eqref{continuity.eq} over all
$i=1,\ldots,N+1$, 
we see that 
\[\rho_{\mathrm{tot}}=\sum_{i=1}^{N+1}
\rho_i\]
is constant in time, and hence $\rho_{\mathrm{tot}}=\sum_{i=1}^{N+1}
(\rho_0)_i=1$. Using this, we can rewrite the implicit formulation of the currents as
\begin{align}\label{flux2}
\nabla \rho_i
=\frac{\rho_i\left(-\sum_{j=1}^NJ_j\right)-
\left(1-\sum_{j=1}^N\rho_j\right)J_i}{D_{i(N+1)}}
+ \sum_{j=1}^N\frac{\rho_iJ_j-\rho_jJ_i}{D_{ij}}
\end{align}
As before, we can define a matrix

\begin{equation}\label{tildeM}
\mathcal M _{ij}(\rho):=\frac{\rho_i}{D_{ij}}-\frac{\rho_i}{D_{i(N+1)}}
-\delta_{ij}\left(\sum_{k=1}^N\frac{\rho_k}{D_{ik}}+\frac{1-\sum_{l=1}^N{\rho_l}}{D_{i(N+1)}}
\right),\qquad i,j=1,\ldots,N.
\end{equation}
From \eqref{ker.im.M}, the matrix $\mathcal M(\rho)$ has full rank, and hence it is invertible. We have
\[J_i= -\sum_{j=1}^NA _{ij} (\rho)\nabla\rho_j\qquad \mbox{with }A(\rho):=-\mathcal M(\rho)^{-1}.\]
\begin{remark}\label{rem3}
The matrix $\mathcal M(\rho)$ is actually independent from the diagonal elements $D_{ii}$.
\end{remark}
\begin{proposition}
Let $s$ be as in~\eqref{eq:definition_of_s}, and let
    $\mathcal M$ be given by \eqref{tildeM}. Then, the
    matrix-valued function $A(\rho):=-\mathcal M(\rho)^{-1}$
    fulfills (H1) and (H2a).
\end{proposition}

\begin{proof}
Let $A(\rho)=-\mathcal M^{-1}(\rho)$. The fact that $\mathcal{M}$ is smooth directly implies that $A$ is
smooth.
Similarly as in the proof of \cite[Lemma~3.2]{JS13}, one can show that
\begin{equation}\label{eq:auxMS}
      \sum_{i=1}^{d}\partial_i w\cdot A(u(w))s''(u(w))^{-1}\partial_i w
      \geq \gamma|\nabla {u(w)}|^2
      \end{equation}
for some $\gamma >0$ and all smooth $w$. 

In order to prove (H2a), we have to show that
\[z\cdot s''(\rho)A(\rho)z \geq \gamma |z|^2\qquad\mbox{for all }z\in \mathbb R^N,\, \rho\in\mathcal D.\]
Let $\rho\in \mathcal D$, ${\bm x}_0\in\Omega$, and $z\in\mathbb R^N$.
We define the following vector-valued function of ${\bm x}$:
        \[w({\bm x}):= s'(\rho)+s''(\rho)z (\bm x-\bm x_0)\cdot \hat e_1,\]
where $\hat e_1$ denotes the unit vector $(1,0,\ldots,0)\in \mathbb
R^d$.
We have 
    \[\partial_i w({\bm x}_0)=\delta_{i1}  s''(\rho)z\]
		and, for $u=(s^\prime)^{-1}$,
		\[
             \partial_i u( w({\bm x}_0))
             = u'(w(\bm x_0))\partial_i w(\bm x_0)
              =u'(w(\bm x_0)) \delta_{i1}s''(\rho) z
                 =u'(w(\bm x_0)) \delta_{i1}s''(u(w(\bm x_0)))z
            =\delta_{i1}z.
          \]
		This, together with~\eqref{eq:auxMS}, implies that 
		\begin{align*}
		z\cdot s''(\rho)A(\rho)z&= 
		(s''(\rho)z)\cdot A(\rho)s''(\rho)^{-1}(s''(\rho)z) \\
		&= 
		\sum_{i=1}^{d}\partial_i w({\bm x}_0)\cdot A(u(w({\bm x}_0)))s''(u(w({\bm x}_0)))^{-1}\partial_i w({\bm x}_0)\\
		&\geq \gamma|\nabla {u(w({\bm x}_0))}|^2 =\gamma |z|^2,
		\end{align*}
		which proves the assertion.
	\end{proof}

For $N=1$, the matrix $\mathcal M(\rho)$ is actually a scalar, which is given by
\begin{align*}
\mathcal M(\rho)=-\frac{\rho_1}{D_{12}}-\frac{1-\rho_1}{D_{12}}=-\frac{1}{D_{12}}.
\end{align*}
Hence, $J_1=D_{12}\nabla\rho_1$. Therefore, in this case the
Maxwell-Stefan system reduces to the heat equation.

	For three species/gases ($N=2$), we have
	\begin{align*}
	\mathcal M(\rho_1,\rho_2)
	&=
	\begin{pmatrix}
	\frac{\rho_1}{D_{11}}-\frac{\rho_1}{D_{13}}-\frac{\rho_1}{D_{11}}-\frac{\rho_2}{D_{12}}-\frac{1-\rho_1-\rho_2}{D_{13}} & \frac{\rho_1}{D_{12}}-\frac{\rho_1}{D_{13}}\\
	\frac{\rho_2}{D_{21}}-\frac{\rho_2}{D_{23}} & 	\frac{\rho_2}{D_{22}}-\frac{\rho_2}{D_{23}}-\frac{\rho_1}{D_{21}}-\frac{\rho_2}{D_{22}}+\frac{1-\rho_1-\rho_2}{D_{23}} 
	\end{pmatrix}
	\\&=-
	\begin{pmatrix}
	\frac1{D_{13}}+\big(\frac{1}{D_{12}}-\frac{1}{D_{13}}\big)\rho_2
	 &\big(\frac{1}{D_{13}}-\frac{1}{D_{12}}\big)\rho_1
	\\
	\big(\frac{1}{D_{23}}-\frac{1}{D_{21}}\big)\rho_2
	&
	\frac1{D_{23}}+\big(\frac{1}{D_{21}}-\frac{1}{D_{23}}\big)\rho_1
	\end{pmatrix}.
	\end{align*}
	Let
	\[d_1:=\frac{1}{D_{13}},\qquad
          d_2:=\frac{1}{D_{23}},\qquad d_3:=\frac{1}{D_{12}},\]
        and recall that $D_{21}=D_{12}$.
	One can verify that
	\[\delta(\rho_1,\rho_2):= \det \mathcal M(\rho_1,\rho_2)= d_1d_2(1-\rho_1-\rho_2)+d_2d_3\rho_1+d_3d_1\rho_2\neq0. \]
    Let $A(\rho)$ denote the inverse of $-\mathcal M(\rho)$. 
    We can rewrite the Maxwell-Stefan equations as the system in
    section~\ref{sec.ThreeGas}.

\subsection{Implicit formulation for the currents}\label{sec:implicit}
In subsection \ref{subsec.maxstef}, we have seen that the Maxwell-Stefan system~\eqref{continuity.eq}-\eqref{flux}, can be written in the form \eqref{eq.classical}, with $f=0$ and $A(\rho)$ being given by the inverse of $-\mathcal M(\rho)$ for 
\begin{equation*}
\mathcal  M _{ij} (\rho):=\frac{\rho_i}{D_{ij}}-\frac{\rho_i}{D_{i(N+1)}}-\delta_{ij}\left(\sum_{k=1}^N\frac{\rho_k}{D_{ik}}+\frac{1-\sum_{l=1}^N{\rho_l}}{D_{i(N+1)}}\right),\qquad i,j=1,\ldots,N.
\end{equation*}
Moreover, we have computed $A(\rho)$ explicitly for $N=1$ and $N=2$. However, for large $N$, it is more complicated to find the explicit formulation for $A(\rho)$. In any case we do not expect a simple formulation in these cases. Therefore, this section provides a space-time Galerkin scheme, which avoids the explicit computation of the inverse of $\mathcal M$.
 
Let $q,p\in\mathbb N$. We consider the following problem:
\begin{multline}\label{MS.galerkin}
\text{Find}\quad w_h^{\epsilon}\in \bm V_h^p,J^\mu\in \bm V_h^q,\ \mu=1,\ldots,d, \quad\text{such that}\\
\begin{aligned}
  0&= \epsilon(\phi^0,w_{h}^{\epsilon})_{{\WeightedSpace}(Q_T)}+
  \int_{\Omega}\phi^0(T)\cdot u(w_h^{\epsilon}(T)) dx
  -\int_{\Omega}\phi^0(0)\cdot\rho_0 dx
-\int_{Q_T}\partial_t \phi^0 \cdot u(w_h^{\epsilon}) dxdt
\\&\qquad-\sum_{\mu=1}^d\left(\int_{Q_T}\partial_{x_\mu} \phi^0 \cdot J^\mu dxdt
+\int_{Q_T}\phi^\mu \cdot\big(\partial_{x_\mu} w_h^{\epsilon}-s''(u(w_h))\mathcal M(u(w_h^{\epsilon}))J^\mu\big) dxdt\right)
\end{aligned}
\\ \forall \phi^0\in \bm V_h^p,\ \phi^\mu\in \bm V_h^q,\ \mu=1,\ldots,d.
\end{multline}

\begin{proposition}\label{Max-Stef-existence}
	Assume that $\rho_0:\Omega\to\overline{\mathcal D}$ is measurable. Then there exists a solution $w_h^\epsilon\in \bm V_h^p ,J^\mu\in \bm V_h^q$, $\mu=1,\ldots,d$ of the method \eqref{MS.galerkin}.
      \end{proposition}

      For the proof of Proposition~\ref{Max-Stef-existence}, we
        need the following lemma.
      
\begin{lemma}\label{Max-Stef-Entropy}
	If $ w_h^{\epsilon}\in \bm V_h^p,J^\mu\in \bm V_h^q,\
        \mu=1,\ldots,d $, solves \eqref{MS.galerkin}, then
	\[
       \epsilon\|w_h^{\epsilon}\|_{{\WeightedSpace}(Q_T)}^2+   \int_{\Omega}s(u(w_h^{\epsilon}(T)))dx +\gamma  \sum_{\mu=1}^d\int_{Q_T}|\mathcal  M(u(w_h^{\epsilon}))J^\mu|^2dxdt
          \leq  \int_{\Omega}s(\rho_0) dx.
        \]
\end{lemma}
\begin{proof}
	We can use $\phi^0=w_h^{\epsilon}$ and $\phi^\mu=0$
        for $\mu=1,\ldots,d$ as  test functions and, similarly to the
        proof of Proposition \ref{prop.nonlin}, we obtain
        that
	\begin{align*}
\epsilon\|w_h^{\epsilon}\|_{{\WeightedSpace}(Q_T)}^2+	\int_{\Omega}s(u(w_h^{\epsilon} (T)))dx - \sum_{\mu=1}^d\int_{Q_T} J^\mu \cdot\partial_{x_\mu} w_h^{\epsilon} dxdt
	\leq  \int_{\Omega}s(\rho_0) dx.
	\end{align*}
	The next step is to use the test functions $\phi^0=0$ and $\phi^\mu=J^\mu$ for $\mu=1,\ldots,d$ to obtain
	\begin{align*}
	\sum_{\mu=1}^d\int_{Q_T}J^\mu\cdot\partial_{x_\mu}
          w_h^{\epsilon}
          dxdt=\sum_{\mu=1}^d\int_{Q_T}J^\mu \cdot s''(u(w_h^{\epsilon}))\mathcal
          M(u(w_h^{\epsilon}))J^\mu dxdt.
 	\end{align*}
	According to assumption (H2a), we know that $s''(v)A(v)$ is positive semi-definite and satisfies
	\[z\cdot s''(v)A(v)z\geq \gamma |z|^2\quad \mbox{for all }z\in\mathbb R^N,\, v\in \mathcal D.\]
	Choosing $v=u(w_h^{\epsilon})$, $z:= \mathcal M(u(w_h^{\epsilon}))J^\mu$, we see that
	\begin{align*}
	\gamma |\mathcal  M(u(w_h^{\epsilon}))J^\mu|^2&\leq J^\mu \cdot \mathcal M(v)s''(v)A(v)\mathcal M(v) J^\mu=- J^\mu \cdot \mathcal M(v)s''(v) J^\mu,
	\end{align*}
        where in the last step we have used that $A(v)$ is the
          inverse of $-\mathcal M(v)$. Thus,
	we conclude that 
    \begin{align*}
\epsilon\|w_h^{\epsilon}\|_{{\WeightedSpace}(Q_T)}^2+\int_{\Omega}s(u(w_h^{\epsilon}(T)))dx +\gamma  \sum_{\mu=1}^d\int_{Q_T}|\mathcal  M(u(w_h^{\epsilon}))J^\mu|^2dxdt
    &\leq  \int_{\Omega}s(\rho_0) dx.
    \end{align*}
\end{proof}

\begin{proof}[Proof of Proposition \ref{Max-Stef-existence}]
	The idea of the proof is to proceed similarly to the proof of
        Proposition \ref{prop.nonlin}. We define the 
        mapping
	\[\Phi:\bm V_h^p\times (\bm V_h^q)^d\to \bm V_h^p\times (\bm V_h^q)^d,\ (v,I^1,\ldots,I^d)\mapsto (w,J^1,\ldots,J^d),\]
	 where $w$ is (uniquely) defined via the equation
\begin{multline*}
  0= \epsilon(\phi^0,w)_{{\WeightedSpace}(Q_T)}+
  \int_{\Omega}\phi^0(T) \cdot u(v(T)) dx
    -\int_{\Omega}\phi^0(0)\cdot\rho^0 dx
-\int_{Q_T}\partial_t \phi^0 \cdot u(v) dxdt
\\-\sum_{\mu=1}^d\int_{Q_T}\partial_{x_\mu} \phi^0\cdot I^\mu dxdt\quad \mbox{for all }\phi^0\in \bm V_h^p,
\end{multline*}
 and $J^\mu$ denotes the unique solution (see below for a justification) of 
 \begin{align}\label{eq.for.J}
\int_{Q_T}\phi^\mu\cdot \partial_{x_\mu} v dxdt=\int_{Q_T}\phi^\mu
   \cdot s''(u(v))\mathcal M(u(v))J^\mu dxdt
\qquad \mbox{for all } \phi^\mu\in \bm V_h^q.
 \end{align}
 Note that the mapping $\Phi$ is well-defined, as \eqref{eq.for.J}
 admits a unique solution for given $v\in \bm V_h^p$ according to the
 Lemma of Lax-Milgram: we see that $\partial_{x_\mu} v\in
 L^2(Q_T)^{N}$ and the matrix ${-s''(u(v))\mathcal M(u(v))\in
   L^\infty(Q_T) ^{N\times N}}$ is positive definite,
 because for all $z\in \mathbb R^N$
 \begin{align*}
   z\cdot \big(-s''(u(v))\mathcal M(u(v))\big)  z
   &=A(u(v))y\cdot s ''(u(v))y
 \\&= y\cdot s''(u(v))A(u(v)) y
 \\&\!\!\stackrel{\text{(H2a)}}\geq \gamma |y|^2= \frac{\gamma}{ \|A(u(v))\|^2}\|A(u(v))\|^2|y|^2
 \\&\geq  \frac{\gamma}{ \|A(u(v))\|^2} |A(u(v))y|^2= \frac{\gamma}{ \|A(u(v))\|^2}|z|^2
 \end{align*}
 for $y:=A(u(v))^{-1}z=-\mathcal M(u(v))z$. Moreover, the mapping
 $\Phi$ is continuous since $A$ and $u$ are continuous. Then by the
 Leray-Schauder fixed-point theorem, we obtain that $\Phi$ admits a
 fixed-point if we can show that the set
	\[\{(w,J^1,\ldots J^d)\in \bm V_h\times (\bm V_h^q)^d: (w,J^1,\ldots J^d)=\sigma\Phi(w,J^1,\ldots J^d), \sigma\in[0,1]\}\]
	is bounded. 	Let $(w,J^1,\ldots,J^d)=\sigma \Phi(w,J^1,\ldots,J^d)$ for $\sigma\in(0,1]$. Similarly to Lemma \ref{Max-Stef-Entropy}, we can prove the entropy estimate
	\[
          \frac{\epsilon}{\sigma}\|w\|_{{\WeightedSpace}(Q_T)}^2+\int_{\Omega}s(u(w(T)))dx +\frac{\gamma}{\sigma}  \sum_{\mu=1}^d\int_{Q_T}|\mathcal  M(u(w))J^\mu|^2dxdt
	\leq  \int_{\Omega}s(\rho_0) dx.\]
	Using that $\sigma\in(0,1]$ is bounded from above yields a uniform bound on $w$ in $\bm V_h^q$ and on $\mathcal M(u(w))J^\mu$ in $L^2(Q_T) ^{N}$. As $\bm V_h^q$ is finite dimensional, we directly obtain that $\|w\|_{L^\infty(Q_T) ^{N}}$ is uniformly bounded. Thus,
	\[\|J^\mu\|_{L^2(Q_T) ^{N}}\leq
          \|A(u(w))\|_{L^\infty(Q_T) ^{N\times N}}\|\mathcal M(u(w))J^\mu\|_{L^2(Q_T) ^{N}}\]
	is also uniformly bounded. As all norms are equivalent on $\bm
        V_h^q$, this directly implies that $J^\mu$ is uniformly bounded in $\bm V_h^q$.	Thus, the Leray-Schauder theorem is applicable and yields that $\Phi$ has a fixed-point, and therefore the scheme \eqref{MS.galerkin} admits a solution.
\end{proof}

\begin{proposition}\label{prop2.MS}
    Let $\rho_0:\Omega\to\overline{\mathcal D}$ be measurable and  $w_h^\epsilon\in \bm V_h^p,J^{\epsilon,\mu}_h\in \bm V_h^q,\ \mu=1,\ldots,d$, be a solution {of} \eqref{MS.galerkin} for $\epsilon, h>0$. 
	Then there exist a solution $\rho$ of \eqref{eq.weak} and sequences $h_i,\epsilon_i\to0$, as $i\to\infty$, such that 
	\[u(w^{\epsilon_i}_{h_i})\to\rho\qquad\mbox{in }L^r(Q_T)\mbox{, as }i\to \infty\]
	for all $r\in[1,\infty)$. Moreover, $\rho$ satisfies the entropy estimate
	\begin{equation}\label{eq.entropy.weak.again}	\int_{\Omega}s(\rho(\tau))dx+\gamma \int_0^\tau\int_{\Omega}\left|\nabla\rho\right|^2 dxdt
	\leq  \int_{\Omega}s(\rho_0) dx
	\end{equation}
	for all $\tau\in(0,T]$, where $|\Omega|$ is the volume of $\Omega$.
\end{proposition}
\begin{proof}
	The proof is analogue to the proof of Proposition
        \ref{prop2}. We only need to replace Proposition \ref{prop3}
        by Lemma~\ref{lemma:last} below.
      \end{proof}
\begin{lemma}[Convergence of the scheme for fixed $\epsilon>0$]\label{lemma:last}
	Let $ w_h\in \bm V_h^p,J^\mu_h\in \bm V_h^q,\ \mu=1,\ldots,d $
        be a solution of \eqref{MS.galerkin}, with fixed $\epsilon>0$.
	Then there exists $\rho\in H^1(Q_T)^N$ with $\rho(t,x)\in\overline{\mathcal D}$ for a.e. $(t,x)\in Q_T$ and $s'(\rho)\in H^1(Q_T)^N$, and a sequence $h_\ell\to0$  such that
	$$\rho_{h_\ell}:= u(w_{h_\ell})\to \rho\qquad\mbox{and}\qquad w_{h_\ell}\to s'(\rho)$$
	strongly in $L^2(Q_T)$ and weakly in ${\WeightedSpace}(Q_T)$. Moreover, $\rho$ solves 
	\eqref{eq.stable} and satisfies the entropy estimate \eqref{eq.entropy-cont}
        for $w=s'(\rho)$.
      \end{lemma}
\begin{proof}
  The fact that $w_{h}$ is uniformly bounded in ${\WeightedSpace}(Q_T)^N$ yields
  that there exists $w\in H^1(Q_T)^N$ and  subsequence ${h_\ell}\to 0$
  such that $w_{h_\ell}\rightharpoonup w$ in ${\WeightedSpace}(Q_T)^N$, due to the
  Banach-Alaoglu theorem, and $w_{h_\ell}\to w$ in $L^2(Q_T)^N$ due to Rellich's theorem. 
	As $u$ is bounded, the dominated convergence theorem entails
        the
        convergence for $\rho_{h_\ell}\equiv u(w_{h_\ell})$
        to $\rho:=u(w)$
        along another subsequence (which we do not relabel).

        For the second part, we note that,
     due to the Banach-Alaoglu theorem and the boundedness of $\mathcal M(u(w_h))J^\mu_h$ in $L^2(Q_T)^{N}$, we know that there exist $\xi^\mu\in L^2(Q_T)^{N}$ such that, for a subsequence (not being relabeled),
	\[\mathcal M(u(w_h))J^\mu_h\rightharpoonup \xi^\mu\qquad\mbox{weakly in }L^2(Q_T)^{N}.\]
	In particular,
	\[J^\mu_h=- A(u(w_h))\mathcal M(u(w_h)){J^\mu_h}\rightharpoonup -A(\rho)\xi^\mu=:J^\mu\qquad\mbox{weakly in }L^r(Q_T)^{N}\]
	for every $r\in[1,2)$.
	Finally, for every $\phi^\mu\in H^1(Q_T)^N$, $j=0,\ldots,d$, there exist $\phi^\mu_{h_\ell}\in \bm V^p_{h_\ell}\cap \bm V^q_{h_\ell}$ such that $\phi^\mu_{h_\ell}\to \phi^\mu$ in ${\WeightedSpace}(Q_T)^N$. 
	Using $\phi^\mu_{h_\ell}$ as a test function in
        \eqref{MS.galerkin},
        in the limit $h_\ell\to0$, we obtain
	 \begin{align*}
           0&= \epsilon(\phi^0,w)_{{\WeightedSpace}(Q_T)}+
              \int_{\Omega}\phi^0(T)\cdot u(w(T))
              dx-\int_{\Omega}\phi^0(0)\cdot \rho_0 dx
	 	-\int_{Q_T}\partial_t \phi^0 \cdot u(w) dxdt
	 	\\&\qquad-\sum_{\mu=1}^d\left(\int_{Q_T}\partial_{x_\mu}
           \phi^0\cdot J^\mu dxdt
	 	+\int_{Q_T}\phi^\mu\cdot \big(\partial_{x_\mu} w-s''(u(w))\mathcal M(u(w))J^\mu\big)dxdt\right),
	 \end{align*}
         as each integral in \eqref{MS.galerkin} converges separately. In particular, by the fundamental lemma of calculus of variations, we see that $\partial_{x_\mu} w=s''(u(w))\mathcal M(u(w))J^\mu$ and equivalently
         \[{J^\mu=\mathcal M(u(w))^{-1}s''(u(w))^{-1}\partial_{x_\mu} w=-A(u(w))u'(w)\partial_{x_\mu} w=-A(u(w))\partial_{x_\mu}u(w),}\]
	  which implies that $\rho$ solves \eqref{eq.stable}. Finally, the entropy inequality is a consequence of Fatou's lemma.
        \end{proof}

\subsection{Numerical Tests}
\label{numerical_tests}

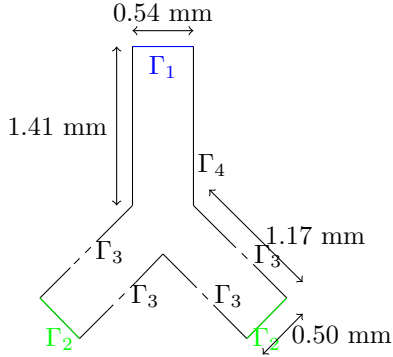
\begin{figure}[!ht]
\centering
\begin{tikzpicture}[scale=1.5]
\def\mypoints{%
    (0,0), (0, 1.41), (-0.54, 1.41), (-0.54, 0.0), (-0.81, -0.27), (-1.09, -0.55), (-1.36, -0.82), (-1.01, -1.18), (-0.76, -0.93), (-0.51, -0.68), (-0.27, -0.43), (-0.02, -0.68), (0.225, -0.93), (0.473, -1.18), (0.827, -0.82), (0.579, -0.57), (0.331, -0.33), (0,0)
};
  \foreach \x [count=\xi, remember=\xi-1 as \xiprev] in \mypoints {
    \node[] (node\xi) at \x {};
    \ifnum\xi>1 %
        \ifthenelse{\xi=13 \OR \xi=10 \OR \xi=6 \OR \xi=17}
        { \draw[dashed] (node\xiprev.center) -- (node\xi.center) node [midway,above,right] {$\Gamma_3$}; }
        { \draw[] (node\xiprev.center) -- (node\xi.center); }
    \fi
  };
\draw[blue] (node2.center) -- (node3.center) node [midway, below] {$\Gamma_1$};
\draw[draw=none] (node2.center) -- (node17.center) node [midway, below] {$\Gamma_4$};
\draw[green] (node14.center) -- (node15.center) node [midway, below] {$\Gamma_2$};
\draw[green] (node7.center) -- (node8.center) node [midway, below] {$\Gamma_2$};

\draw [<->] ([yshift=4 pt]node2.center) -- ([yshift=4 pt]node3.center) node [midway,above] {0.54~mm};
\draw [<->] ([xshift=-4 pt]node3.center) -- ([xshift=-4 pt]node4.center) node [midway,left] {1.41~mm};
\draw [<->] ([xshift=4 pt,yshift=-4 pt]node14.center) -- ([xshift=4 pt,yshift=-4pt]node15.center) node [midway,right] {0.50~mm};
\draw [<->] ([xshift=4 pt,yshift=4 pt]node15.center) -- ([xshift=4 pt,yshift=4pt]node18.center) node [midway,right] {1.17~mm};
\end{tikzpicture}
\caption{Computational domain for the lung model.}
\label{fig.lungdom}
\end{figure}

We again turn to \cite[Sec. 3]{boudin} for numerical results we can compare our method to.
This time, we consider a model for the lung.
The computational domain resembles {a} branch of the tree structure found in the bottom of the lung. 
The domain, depicted in \Cref{fig.lungdom}, consists of the inflow, $\Gamma_1$, on top, the outflow, $\Gamma_2$,
located on the bottom of the two branches, and the alveoli, $\Gamma_3$, located in the middle of each of the branches.
The remaining boundary $\Gamma_4$ is a wall where nothing goes in or out. 
Opposed to the domain presented in the reference, we consider the branches of the lung to be symmetrical and perpendicular to each other. The paper does not mention the angle between the branches used there.
Also the size of the alveoli is left unspecified in the paper. Here, we split the boundary of the branches into three equal parts, with the alveoli ($\Gamma_3$) in the middle.
On $\Gamma_1,\Gamma_2,\Gamma_3$ we impose Dirichlet boundary conditions to model the gas exchange with the other parts of the lung. 
On the wall, $\Gamma_4$, we take homogeneous Neumann boundary conditions.


We make use of the implicit formulation \eqref{MS.galerkin} to find the numerical solution. 
To incorporate the Dirichlet boundary condition, we use Nitsche's
  method and add to \eqref{MS.galerkin} the following terms:
\begin{align*}
\sum_{\mu=1}^d\int_{(0,T)\times\Gamma_D} J^\mu \nu^\mu \cdot \phi^0
+\int_{(0,T)\times\Gamma_D} (u(w)-\rho_D)\cdot\phi^\mu \nu^\mu 
+\int_{(0,T)\times\Gamma_D} \eta {h_s}^{-1}  (u(w)-\rho_D) \cdot \phi^0
\end{align*}
for a parameter $\eta>0$, $h_s$ being the spatial mesh size, on the Dirichlet boundary $\Gamma_D$. In the examples below, we use $\eta=1$.
The first term comes from the integration by parts. 
The second and third terms are productive zeros that weakly enforce the Dirichlet boundary condition, and are chosen such that they agree with Nitsche's method for the heat equation in the degenerative case.

\subsubsection{Diffusion of air}
In the following example, compare \cite[Sec. 3.4]{boudin}, we choose alveolar air as initial condition and as the Dirichlet data on the outflow and alveoli.
On the inflow boundary we choose humidified air as Dirichlet data.
See \Cref{tab.lunggas} for the gas components of the different types of air, and \Cref{tab.lungdiff} for the diffusion coefficients.

\begin{table}[!ht]
\centering
\caption{Components of the different gas mixtures.}
\label{tab.lunggas}
\begin{tabular}{ l | c c c}
 & Humidified air & Alveolar air & Alveolar heliox\\ \hline
    Nitrogen & 0.7409 & 0.7490 & 0.0000\\ 
    Oxygen & 0.1967 & 0.1360 & 0.1360\\  
    Carbon dioxide & 0.0004 & 0.0530 & 0.0530\\  
    Water & 0.0620 & 0.0620 & 0.0620 \\  
    Helium & 0.0000 & 0.0000 & 0.7490   
\end{tabular}
\end{table}

\begin{table}[!ht]
\centering
\caption{Diffusion coefficients of the different gases.}
\label{tab.lungdiff}
\begin{tabular}{ l | c c c c}
 & Oxygen & Carbon dioxide & Water & Helium\\ \hline
    Nitrogen & 21.87 & 16.63 & 23.15 & 74.07\\ 
    Oxygen & & 16.40 & 22.85 & 79.07\\  
    Carbon dioxide &  & & 16.02 & 63.45 \\
    Water & & & & 90.59
\end{tabular}
\end{table}

\begin{figure}[!ht]  
\centering
\resizebox{0.4\linewidth}{!}{
    \begin{tikzpicture}
    \begin{axis}
    [ xlabel={Time (s)},
    ylabel={Oxygen},
    ymajorgrids=true,
    grid style=dashed,
    cycle list name=paulcolors,
    yticklabel style={
        /pgf/number format/precision=3,
    },
    xticklabel style = {
        /pgf/number format/fixed,
        /pgf/number format/precision = 3
    }
    ]
    \addplot table [mark=none, x=t, y=O2, col sep=comma] {./lung2.csv};
    \end{axis}
    \end{tikzpicture}
}
\resizebox{0.4\linewidth}{!}{
    \begin{tikzpicture}
    \begin{axis}
    [ xlabel={Time (s)},
    ylabel={Carbon dioxide},
    ymajorgrids=true,
    grid style=dashed,
    cycle list name=paulcolors,
    xticklabel style = {
        /pgf/number format/fixed,
        /pgf/number format/precision = 3
    }
    ]
    \addplot table [mark=none, x=t, y=CO2, col sep=comma] {./lung2.csv};
    \end{axis}
    \end{tikzpicture}
}
\caption{Numerical results of the mole fractions Oxygen and Carbon dioxide inside the lung for air mixture.}
\label{fig.lungresults}
\end{figure}
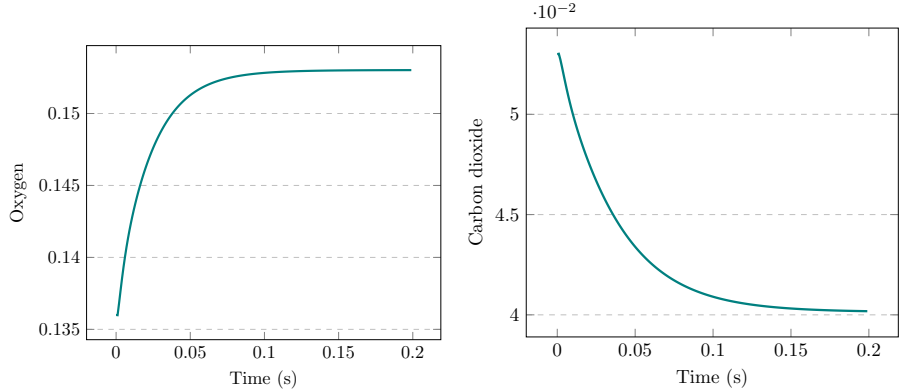

Since there is no helium present we can reduce the number of species involved, setting $N=3$. 
For the numerical calculations we choose spatial mesh size $h_s=0.3$ and measure the value of the gas every 0.001 seconds. 
The discrete system is not ill-conditioned and we are able to choose $\epsilon=0$.
In \Cref{fig.lungresults} we show the numerical results for Oxygen and Carbon dioxide as the other gases stay (almost) constant.
Both converge to their equilibrium value. 
Comparing the results to \cite{boudin}, we can see that the
  equilibrium value slightly differs, which is likely due
  to the symmetry of the domain
  and size of the alveoli. 

\subsubsection{Diffusion of air/heliox}
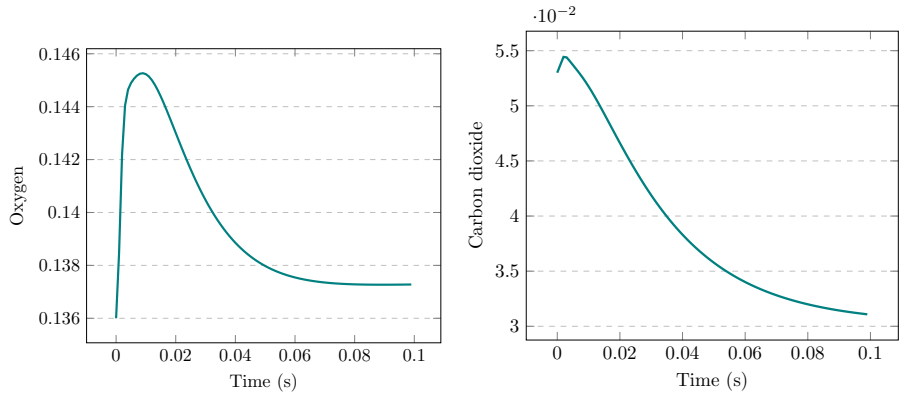
\begin{figure}[!ht]  
\centering
\resizebox{0.4\linewidth}{!}{
    \begin{tikzpicture}
    \begin{axis}
    [ xlabel={Time (s)},
    ylabel={Oxygen},
    ymajorgrids=true,
    cycle list name=paulcolors,
    grid style=dashed,
    yticklabel style={
        /pgf/number format/precision=3,
    },
    xticklabel style = {
        /pgf/number format/fixed,
        /pgf/number format/precision = 3 
    }
    ]
    \addplot table [mark=none, x=t, y=O2, col sep=comma] {./lung3.csv};
    \end{axis}
    \end{tikzpicture}
}
\resizebox{0.4\linewidth}{!}{
    \begin{tikzpicture}
    \begin{axis}
    [ xlabel={Time (s)},
    ylabel={Carbon dioxide},
    ymajorgrids=true,
    grid style=dashed,
    cycle list name=paulcolors,
    xticklabel style = {
        /pgf/number format/fixed,
        /pgf/number format/precision = 3
    }
    ]
    \addplot table [mark=none, x=t, y=CO2, col sep=comma] {./lung3.csv};
    \end{axis}
    \end{tikzpicture}
}
\caption{Numerical results of the mole fractions Oxygen and Carbon dioxide inside the lung for air/heliox mixture.}
\label{fig.lungresults2}
\end{figure}

Next, we try to reproduce the results form \cite[Sec. 3.5]{boudin}.
We consider alveolar heliox as initial condition. 
As the Dirichlet data on the outflow and alveoli, we also choose alveolar heliox, whereas we put humidified air on the inflow.
The discrete system is very ill-conditioned due to the gas components taking zero values. 
In order for the solver to converge, we had to choose {$\epsilon=10^{-3}$}.
Furthermore, to avoid the singularity of the entropy density, we adjust the helium content in air and the nitrogen content in heliox to be $10^{-7}$, subtracting the same amount of water, in order to keep them summing to one. 
Note that this is not unreasonable, for example, the correct amount of helium in air is about $5.3\cdot 10^{-7}$.
With these adjustments, the solver converges.
The numerical results are shown in \Cref{fig.lungresults2}.
Both oxygen and carbon dioxide levels rise above the values in provided gas mixtures, before they start to decrease towards the equilibrium value. 
This is the expected behavior.
However, the maximum values reached here are slightly lower
than the ones found in \cite{boudin}. 
This can be attributed to the perturbations of the zero
  concentrations and, as already seen, to the approximation of the geometry.

\section{Conclusions}\label{sec:conclusions}
We have presented and analyzed a continuous space-time Galerkin method
for cross-diffusion systems in entropy variable formulation, proving existence and convergence of
discrete solutions, as well as existence of a weak solution of the
continuous problem using the space-time approach.
{As opposed to time-stepping schemes, this approach provides an easy
way to increase the approximation order simultaneously in space and time,
makes $hp$-refinement in space-time possible, without the need for a
globally fixed time-step size, and delivers numerical 
solutions, which can be evaluated at arbitrary points in time.
Furthermore, at the same time, positivity and boundedness of the solutions are preserved also at the discrete level.}

In the numerical examples, we have observed optimal convergence rates, given that the solution stays away from the singularities of the entropy. 
Lifting this restriction could be the topic of future works. 
{Also, a more efficient numerical treatment of the space-time system is of interest, for example using a suitable preconditioner, a fine tuned solver, and making use of the mesh structure when using a tensor-product mesh.}

\section*{Acknowledgments}
    All authors have been supported by the Austrian Science Fund
    (FWF) through the project F~65.
    I. Perugia and P. Stocker have also been supported by the FWF
    through the projects P~29197-N32 and~W1245, respectively.

\bibliography{bib}{}
\bibliographystyle{siamplain}
\end{document}

%% file: pgfslopetriangle.tex
\newcommand{\logLogSlopeTriangle}[5]
{

    \pgfplotsextra
    {
        \pgfkeysgetvalue{/pgfplots/xmin}{\xmin}
        \pgfkeysgetvalue{/pgfplots/xmax}{\xmax}
        \pgfkeysgetvalue{/pgfplots/ymin}{\ymin}
        \pgfkeysgetvalue{/pgfplots/ymax}{\ymax}

        \pgfmathsetmacro{\xArel}{#1}
        \pgfmathsetmacro{\yArel}{#3}
        \pgfmathsetmacro{\xBrel}{#1-#2}
        \pgfmathsetmacro{\yBrel}{\yArel}
        \pgfmathsetmacro{\xCrel}{\xBrel}

        \pgfmathsetmacro{\lnxB}{\xmin*(1-(#1-#2))+\xmax*(#1-#2)} 
        \pgfmathsetmacro{\lnxA}{\xmin*(1-#1)+\xmax*#1} 
        \pgfmathsetmacro{\lnyA}{\ymin*(1-#3)+\ymax*#3} 
        \pgfmathsetmacro{\lnyC}{\lnyA+#4*(\lnxA-\lnxB)}
        \pgfmathsetmacro{\yCrel}{\lnyC-\ymin)/(\ymax-\ymin)} 

        \coordinate (A) at (rel axis cs:\xArel,\yArel);
        \coordinate (B) at (rel axis cs:\xBrel,\yBrel);
        \coordinate (C) at (rel axis cs:\xCrel,\yCrel);

        \draw[#5]   (A)-- 
                    (B)-- node[pos=1.2,left=-3pt] {#4}
                    (C)-- 
                    cycle;
    }
}